\documentclass{article}
\usepackage{indentfirst}
\usepackage{arxiv}
\usepackage[utf8]{inputenc} 
\usepackage[T1]{fontenc}    
\usepackage{booktabs}       
\usepackage{amsfonts}       
\usepackage{nicefrac}       
\usepackage{microtype}      
\usepackage{lipsum}		
\usepackage{graphicx}
\usepackage[numbers]{natbib}
\usepackage{hyperref}
\usepackage[normalem]{ulem}
\hypersetup{colorlinks}
\usepackage{amsmath}
\usepackage{amsthm}
\usepackage{float}
\restylefloat{table}
\usepackage{scalerel}
\usepackage[ruled,vlined,linesnumbered,noresetcount,algosection]{algorithm2e}
\numberwithin{equation}{section}
\usepackage{chngcntr}
\counterwithin{figure}{section}
\usepackage{mathtools}
\usepackage{caption}
\usepackage[title]{appendix}
\DeclarePairedDelimiter\ceil{\lceil}{\rceil}

\title{ACTIVE LEARNING BASED SAMPLING FOR\\ HIGH-DIMENSIONAL NONLINEAR PARTIAL DIFFERENTIAL EQUATIONS}

\author{{Wenhan Gao}\\
Department of Applied Mathematics and Statistics, Stony Brook University, Stony Brook, NY 11794, USA
\\
	\texttt{wenhan.gao@stonybrook.edu} \\
\AND
    {Chunmei Wang}\\
   	Department of Mathematics, University of Florida, Gainesville, FL 32611, USA \\
	\texttt{chunmei.wang@ufl.edu}\\
}

\date{}
\renewcommand{\headeright}{}
\renewcommand{\undertitle}{}
\renewcommand{\shorttitle}{ACTIVE LEARNING BASED SAMPLING }

\newtheorem{thm}{Theorem}[section]

\newtheorem{definition}{Definition}[section]

\begin{document}
\maketitle

\begin{abstract}
	The deep-learning-based least squares method has shown successful results in solving high-dimensional and non-linear partial differential equations (PDEs). However, this method usually converges slowly. To speed up the convergence  of this approach, an active-learning-based sampling algorithm is proposed in this paper. This algorithm actively chooses the most informative training samples from a probability density function based on residual errors to facilitate error reduction. In particular, points with larger residual errors will have more chances of being selected for training. This algorithm imitates the human learning process: learners are likely to spend more time repeatedly studying mistakes than 
other tasks they have correctly finished. A series of numerical results are illustrated to demonstrate the effectiveness of our active-learning-based sampling in high dimensions to speed up the convergence of the deep-learning-based least squares method.
\end{abstract}


\subjects{65M75 \and 35J25 \and 65N99}

\section{Introduction}
\subsection{Problem Statement}
High-dimensional partial differential equations are widely used in modeling complicated phenomena; e.g., the Hamilton-Jacobi-Bellman (HJB) equation  \cite{HJB} in control theory, the Schr$\ddot{o}$dinger equation  \cite{sch_equation} in quantum mechanics. There are various  traditional numerical methods for solving PDEs, especially for low-dimensional and linear problems. However, when it comes to solving high-dimensional problems, the curse of dimensionality becomes a major computational challenge for many traditional methods such as  the finite difference method  \cite{FDM} and the finite element methods where the computational complexity is exponential to the number of dimensions. With that being said, traditional methods oftentimes are computationally intractable in solving high-dimensional problems.

In approximation theory, deep neural networks (DNNs) have been  a more effective tool for function approximation than traditional approximation tools such as finite elements and wavelets. In \cite{yarotsky1,yarotsky2,shijun2,shijun3,shijun6,hon2021simultaneous}, it has been shown that ReLU-type DNNs can provide more attractive approximation rates than traditional tools to approximate continuous functions and smooth functions with explicit error characterization \cite{shijun2,shijun3,shijun6} in the $L^\infty$-norm \cite{yarotsky1,yarotsky2,shijun2,shijun3,shijun6} and in the $W^{s,p}$-norm for $p\in[1,\infty)$ \cite{hon2021simultaneous}. If target functions lie in the space of functions with special integral representations \cite{barron1993,weinan1,weinan2,siegel2021sharp,SIEGEL2020313,JCM-39-801}, ReLU-type DNNs can lessen the curse of dimensionality in function approximation. Armed with advanced activation functions, DNNs are able to conquer the curse of dimensionality in approximation theory   for approximating continuous functions \cite{shijun4,shijun5,yarotsky3,shijun7}. In generalization theory, it was shown that DNNs can achieve a dimension-independent error rate for solving PDEs \cite{doi:10.1137/19M125649X,LuoYang2020,lu2021priori,chen2021representation}. These theoretical results have justified the application of DNNs to solve high-dimensional PDEs recently in \cite{Han8505,raissi2017physicsI, raissi2017physicsII, Deep_Ritz, DLPDE1, DLPDE2l, Deep_Ritz, IntNet, DLPDE3}. Two main advantages of deep-learning-based methods presented in these studies can be summarized as follows: firstly, the curse of dimensionality can be weakened or even be overcome in certain classes of PDEs \cite{arnulf1,arnulf2}; secondly, deep-learning-based PDE solvers are mesh-free without tedious mesh generation for complex domains in traditional solvers. Thus, deep-learning-based methods have shown tremendous potential to surpass other numerical methods especially in solving high-dimensional PDEs in complex domains.

Deep-learning-based PDE solvers have been applied to solve problems in many areas of practical engineering and life sciences   \cite{other1, other2, other3}. However, deep-learning-based solvers usually are computationally expensive meaning that the convergence speed is usually slow. To reach a desired precision, the training process can take a very long time. Therefore, it is imperative to speed up the convergence. A recent study  \cite{gu2021selectnet} proposed a self-paced learning framework, that was inspired by an algorithm named "SelectNet"  \cite{liu2019selectnet} originally  for image classification, to ease the issue of slow convergence. In   \cite{gu2021selectnet}, new neural networks called selection networks were introduced to be trained simultaneously with the residual model based solution network. The selection network learned to provide weights to training points to guide the solution network to learn more from points that have larger residual errors.

In this paper, we aim to address the aforementioned issue of slow convergence by introducing active-learning-based adaptive sampling techniques to sample the most informative training examples for networks to be trained on. This algorithm is motivated by the active learning  approach primarily used in supervised learning in computer vision  \cite{AL1, AL2, AL3, AL4},  which chooses high-quality data to label so that the network learns faster and less human labeling labor is encountered.

  \subsection{Related work}
\subsubsection{Relevant Literature Review}\label{sec:relevant lit.}
\label{sec:2}
Mesh-based neural networks are proposed to approximate solution operators for a specific type of differential equations \cite{LEE1990110, Lagaris_1998, MALEK2006260, meshbased_1, meshbased_2}. However, these methods become computationally intractable in high dimensions since the degrees of freedom will explode exponentially as the number of dimensions increases. Therefore, the mesh-free model   \cite{deep_XDE,meshfree2,E_2017, DLPDE2l,raissi2017physicsI,meshfree1,WAN,chen2021friedrichs} is the key to solving high-dimensional PDEs. Mesh-free models use sampling of points and the Monte-Carlo method to approximate the objective loss function to avoid generating meshes that are computationally restrictive in high dimensions. In an early study \cite{DLPDE4}, PDEs are solved via deep neural networks by directly minimizing the degree to which the network approximation violates the PDE at prescribed collocation points, and boundary conditions are treated as a hard constraint by constructing boundary governing functions that satisfy the boundary conditions although designing such functions usually is difficult in complicated PDE problems. In a later study concerning boundary value problems  \cite{DLPDE5}, boundary residual errors are also taken into account in a single loss function to enforce boundary conditions. The latter is more commonly used since it does not require auxiliary functions that are problem-dependent and difficult to construct. 
 
 Active learning methods are not commonly utilized in deep-learning-based PDE solvers. In  \cite{activePDE1}, active learning is used to select ``additional observation locations" to "minimize the variance of the PDE coefficient" and to "minimize the variance of the solution of the PDE". In  \cite{deep_XDE}, the residual-based adaptive refinement (RAR) is proposed to adaptively select large residual points. RAR first chooses uniform points and acquires residual errors at these points, and then replicates a certain number of large residual points to the set for training. RAR and our active-learning-based adaptive sampling (to be introduced in Section \ref{sec:3.2}) share the same methodology of selecting large residual error points based upon uniform points. The difference is that our adaptive sampling ranks all points. Low residual points will less likely be selected whereas RAR still keeps them. Large residual error points are also ranked, i.e., even if two points both have large residuals, the one with larger residual error is more likely to be selected. In contrast, RAR does not differentiate them and chooses the largest $k$ residual points. In \cite{Adaptive1}, various adaptive strategies have also been introduced. For example, in their adaptive-sampling-in-time approach, collocation points initially come from a small specified time interval and this interval gradually expands. As shown in \cite{Adaptive1}, this approach considerably improves the convergence and overall accuracy. However, the number of training points accumulates in this approach and the computation becomes more expensive. In our adaptive sampling approach, the number of training points does not grow; therefore, in some cases, our adaptive sampling is a better fit although the approaches in \cite{Adaptive1} are already highly sophisticated.

\subsubsection{Deep-Learning-based Least Squares Method}
\label{sec:2.1}

\noindent Consider the following boundary value problem for simplicity to introduce the deep-learning-based least squares method: find the unknown solution $u(\boldsymbol{\boldsymbol{x}})$ such that 
\begin{equation}
\label{eq:2.1}
\left\{\begin{array}{ll}
\mathcal{D}u(\boldsymbol{\boldsymbol{x}})=f(\boldsymbol{\boldsymbol{x}}), & \text { in } \Omega, \\
\mathcal{B}u(\boldsymbol{\boldsymbol{x}})=g(\boldsymbol{\boldsymbol{x}}), & \text { on } \partial \Omega,
\end{array}\right.
\end{equation}
 where $\partial \Omega$ is the boundary of the domain, $\mathcal{D}$ and $\mathcal{B}$  are   differential operators in $\Omega$ and on $\partial \Omega$, respectively. The goal is to train a neural network, denoted by $\phi(\boldsymbol{x}; \boldsymbol{\theta})$, where $\boldsymbol{\theta}$ is the set of network parameters, to approximate the ground truth solution  $u(\boldsymbol{x})$ of the PDE \eqref{eq:2.1}.
In the least squares method (LSM), the PDE is solved by finding the optimal set of parameters $\boldsymbol{\boldsymbol{\theta}}^\ast$ that minimizes the loss function; i.e.,
\begin{equation}
\label{eq:2.3}
\begin{aligned}
\boldsymbol{\theta}^\ast = \underset{\boldsymbol{\theta}}{\arg \min}\text{ }\mathcal{L}(\boldsymbol{\theta}) 
:&=\underset{\boldsymbol{\theta}}{\arg \min}\text{ } \|\mathcal{D} \phi(\boldsymbol{x} ; \boldsymbol{\theta})-f(\boldsymbol{x})\|_{2}^{2}+\lambda\|\mathcal{B} \phi(\boldsymbol{x} ; \boldsymbol{\theta})-g(\boldsymbol{x})\|_{2}^{2}\\
&=\underset{\boldsymbol{\theta}}{\arg \min}\text{ } \mathbb{E}_{\boldsymbol{x} \in \Omega}\left[|\mathcal{D} \phi(\boldsymbol{x} ; \boldsymbol{\theta})-f(\boldsymbol{x})|^{2}\right]+\lambda \mathbb{E}_{\boldsymbol{x} \in \partial \Omega}\left[|\mathcal{B} \phi(\boldsymbol{x} ; \boldsymbol{\theta})-g(\boldsymbol{x})|^{2}\right]\\
&\approx \underset{\boldsymbol{\theta}}{\arg \min}\text{ } \frac{1}{N_{1}} \sum_{i=1}^{N_{1}}\left|\mathcal{D} \phi\left(\boldsymbol{x}_{i} ; \boldsymbol{\theta}\right)-f\left(\boldsymbol{x}_{i}\right)\right|^{2}+ \frac{\lambda}{N_{2}} \sum_{j=1}^{N_{2}}\left|\mathcal{B} \phi\left(\boldsymbol{x}_{j} ; \boldsymbol{\theta}\right)-g\left(\boldsymbol{x}_{j}\right)\right|^{2},
\end{aligned}
\end{equation} 
where $\lambda$ is a positive hyper-parameter that weights the boundary loss. The last step is a Monte-Carlo approximation with $\boldsymbol{x}_i \in \Omega$ and $\boldsymbol{x}_j \in \partial \Omega$ being $N_1$ and $N_2$ allocation points sampled from the respective probability densities that $\boldsymbol{x}_i$ and $\boldsymbol{x}_j$ follow. This model will be referred to as the basic model in this paper. For a time-dependent PDE, the temporal coordinate can be regarded as another spatial coordinate to build a similar least squares model.

\subsection{Organization}
This paper is structured as follows.
The active-learning-based adaptive sampling will be introduced in Section \ref{sec:3}. The detailed numerical implementation of the proposed sampling method is discussed in  Section \ref{sec:4}. In Section \ref{sec:5}, extensive numerical experiments will be provided to demonstrate the efficiency of the proposed method.

\section{Active-Learning-based Adaptive Sampling}
\label{sec:3}
\subsection{Overview of Active Learning} \label{sec:overview_of_AL}
Active learning \cite{AL5, AL6} is a machine learning method in which the learning algorithm inspects unlabeled data and interactively chooses the most informative data points to learn. Active learning methodology has primarily been applied in computer vision tasks  \cite{AL1}, and has shown an extraordinary performance. Active learning for deep learning proceeds in rounds. In each iteration/round, the current model serves to assess the informativeness of training examples. It is usually used for supervised or semi-supervised learning in which there is a large amount of unlabeled data and the algorithm chooses the most informative data points to get labeled for training. There are two most common metrics to measure informativeness: uncertainty \cite{uncertainty} and diversity \cite{diversity}. In the uncertainty sampling, the algorithm tries to find training examples that the current model is most uncertain about, as a proxy of the model output being incorrect. For example, in a classification problem, this uncertainty can be measured by the current model prediction; if for a training example, the model output probability distribution over classes has similar values as opposed to high probability for a single class, then the current model is uncertain about this particular training example. In the diversity sampling  \cite{diversity2}, the model seeks to identify training examples that most represent the diversity in the data.  It is worth mentioning that query by committee (QBC) \cite{QBC1} is also a common active learning approach. Our proposed algorithm will be a diversity sampling approach that focuses on training examples that the neural network model is more uncertain about. It should be pointed out that diversity sampling and uncertainty sampling are not contradictory to each other. One can combine these two approaches to Active Learning.

\subsection{Adaptive Sampling (AS)}\label{sec:3.2}
In supervised or semi-supervised learning, the uncertainty sampling algorithm chooses data points that the current model is most uncertain about to get labeled and to be trained on. The least squares method for solving PDEs is an unsupervised learning approach because there is no label in the loss function. However, the idea of uncertainty sampling can still be applied to choose the most informative allocation points to learn. Hence, inspired by the uncertainty sampling, we propose the active-learning-based adaptive sampling in this paper to select high-quality train examples to speed up the convergence of PDE solvers with the basic least squares idea in \eqref{eq:2.3}.

In our adaptive sampling, the objective is to preferentially choose allocation points with larger absolute residual errors. Intuitively, one can think of large residual error at a point as a proxy of the model being wrong to a greater extent at this particular point. The fundamental methodology of adaptive sampling is to choose from a biased importance distribution that attaches higher priority to important volumes/regions of the domain. For simplicity,   the absolute residual error of the network approximation at a point $\boldsymbol{x}$ is defined as follows:  
\begin{equation}\label{eq:3.1}
\mathcal{R}_{abs}(\boldsymbol{x}) = 
\left\{\begin{array}{ll}
|\mathcal{D}\phi(\boldsymbol{x}) - f(\boldsymbol{x})| & \text { if } \boldsymbol{x} \in \Omega, \\
|\mathcal{B}\phi(\boldsymbol{x}) - g(\boldsymbol{x})| & \text { if } \boldsymbol{x} \in \partial \Omega,
\end{array}\right.
\end{equation}
 where $\phi$ is the solution network. The adaptive sampling is thus proposed to choose allocation points following a distribution, denoted by $\pi$, 
 given by the probability density function   proportional to the distribution of residual errors, that is $q(\boldsymbol{x}) \propto \mathcal{R}^p_{abs}(\boldsymbol{x})$ in $\Omega$ and on $\partial \Omega$, respectively. The following generalized density function can be applied to choose allocations in $\Omega$ or on $\partial \Omega$, respectively. More precisely, we define
\begin{equation}\label{eq:3.2}
q(\boldsymbol{x}) = \frac{\mathcal{R}^p_{abs}(\boldsymbol{x})}{NC}
\end{equation}
 where $NC = \int_{\Omega \text{ or } \partial \Omega} \mathcal{R}^p_{abs}(\boldsymbol{x}) d\boldsymbol{x}$ is an unknown normalizing constant,  $p$ is a non-negative constant exponent that controls the effect of adaptive sampling. With larger $p$, larger residual error points are more likely to be sampled; when $p = 0$, there is no effect of adaptive sampling. 

As one may have noticed, the unknown normalizing constant is difficult to calculate. Therefore, one may not be able to directly sample points from $q(\boldsymbol{x})$. Two techniques to simulate  observations from $q(\boldsymbol{x})$ are reviewed in the following.

\subsubsection{Metropolis Hastings Sampling}
The Metropolis-Hastings (MH)  algorithm  \cite{MH1} is  a  commonly used Markov  chain  Monte  Carlo (MCMC) method \cite{MCMC1}.   Markov chain   comes  in  because  the  next sample  depends  on  the  current sample.  Monte Carlo comes in because it  simulates observations from a target distribution that is difficult to directly sample from. The Metropolis-Hastings algorithm will be adopted to select points from the target distribution $q(\boldsymbol{x})$.
\begin{algorithm}
\caption{Metropolis Hastings Sampling}\label{alg:3.1}
\SetAlgoLined
\setcounter{AlgoLine}{0}
\SetKwInput{KwResult}{Result}
\KwResult{$N_1$ points in $\Omega$ for training}
\SetKwInput{KwRequire}{Require}
\KwRequire{PDE (\ref{eq:2.1}); the current solution net $\phi(\boldsymbol{\theta})$  }
 Initialize a random point $x_0 \in \Omega$, $j := 0$ and $m = N_1 + b$, where b is a positive integer \;
 \While{j < m}{
        $j:= j+1$ \;
        Draw $x_{candidate} \sim Unif(\Omega)$, and $u \sim Unif(0, 1)$ \;
        acceptance$\_$rate = $\frac{\mathcal{R}^p_{abs}(x_{candidate})}{\mathcal{R}^p_{abs}(x_{l-1})}$ \;
        \uIf{$u <$ acceptance\_rate}{
        $x_j = x_{candidate}$ \;
        }
        \Else{
        $x_{j} = x_{j-1}$  \;
        }
  }
\Return last $N_1$ points\;
\end{algorithm}
\noindent Note that all algorithms sample  points in $\Omega$ for training, and
the same logic follows for sampling on $\partial \Omega$.

In Algorithm  \ref{alg:3.1},    the first $b$ points will be discarded which are called "burn-ins" where $b$ is a non-negative integer. It is a computational trick to correct the bias of early samples. In practice, $b$ can be $0$. However, there is no rigorous justification for burn-ins based on the best of our knowledge. Note that in actual implementation,  Algorithm  \ref{alg:Metropolis Hastings Sampling}  in a vectorized version following the coding style in Python will be deployed.

\begin{algorithm}
\caption{Vectorized Metropolis Hastings Sampling}\label{alg:Metropolis Hastings Sampling}
\SetAlgoLined
\setcounter{AlgoLine}{0}
\SetKwInput{KwResult}{Result}
\KwResult{$N_1$ points in $\Omega$ for training}
\SetKwInput{KwRequire}{Require}
\KwRequire{PDE (\ref{eq:2.1}); the current solution net $\phi(\boldsymbol{x};\boldsymbol{\theta})$ }
 Generate an array of $N_1 + b$ uniform points $ \boldsymbol{X}\subset\Omega$ \;
 RE$\_$array := $\mathcal{R}^p_{abs}(\boldsymbol{X}) = |\mathcal{D}\phi(\boldsymbol{X};\boldsymbol{\theta}) - f(\boldsymbol{X})|^p$ \;
 Generate an array of $N_1 + b$ uniform points $U$ following a uniform distribution on $(0,1)$  \;
 
 \For{$i$ in range $(0 , N_1+b-1)$ } {
        \If{RE$\_$array[i+1] < RE$\_$array[i]}{
        acceptance$\_$rate = $\frac{\text{RE$\_$array[i+1]}}{\text{RE$\_$array[i]}}$ \;
            \If{acceptance$\_$rate < U[i+1]}{
            $\boldsymbol{X}[i+1] =  \boldsymbol{X}[i]$
            }
        }
 
  }
\Return the last $N_1$ points in $\boldsymbol{X}$\;
\end{algorithm}

\subsubsection{Self-normalized Sampling}
The Metropolis-Hastings Algorithm sometimes can be computationally expensive, and the convergence of this algorithm can be slow. Therefore, an alternative algorithm is proposed to sample training examples. This algorithm is faster than the MH algorithm because this algorithm is entirely parallel in generating points as opposed to the MH algorithm in which each sample point is based on the previous sample; therefore, the MH algorithm is not parallel. Our proposed new sampling technique will be used in the numerical experiment. 

\begin{algorithm}
\caption{Self-normalized Sampling}\label{alg:Self-normalized Sampling}
\SetAlgoLined
\setcounter{AlgoLine}{0}
\SetKwInput{KwResult}{Result}
\KwResult{$N_1$ points in $\Omega$ for training}
\SetKwInput{KwRequire}{Require}
\KwRequire{PDE (\ref{eq:2.1}); the current solution net $\phi(\boldsymbol{x};\boldsymbol{\theta})$}
 Generate $N_1$ uniformly distributed points $\big\{\boldsymbol{x}_i\big\}_{i = 1}^{N_1} \subset \Omega$; denote by $\boldsymbol{X}$\;
 RE$\_$array := $\mathcal{R}^p_{abs}(\boldsymbol{X}) = |\mathcal{D}\phi(\boldsymbol{X};\boldsymbol{\theta}) - f(\boldsymbol{X})|^p$ \;
 SRE = sum(RE$\_$array) \;
 probability$\_$array =  $\frac{\text{ RE$\_$array}}{\text{SRE}}$ \;
 Generate $N_1$ points, $\boldsymbol{X}\_$training, following the discrete probability density function $f($RE$\_$array[$i$]$)$ =  probability$\_$array[$i$] \;
 \Return $\boldsymbol{X}\_$training \;
\end{algorithm}

This algorithm is named "Self-normalized Sampling"
because the discrete probability density function obtained by this algorithm normalizes itself. It is well-known that the marginal distribution of points generated by the Metropolis-Hastings algorithm will converge to the target distribution; the self-normalized sampling algorithm does efficiently generate a set of points that focuses on large-residual-error volumes in the experiment. Moreover, in Appendix \ref{self-norm}, we will discuss why this self-normalized algorithm makes sense.

In general, self-normalized sampling is less computationally expensive than the Metropolis-Hastings  algorithm. However, in some rare cases, the residual error is only large in some extremely small volumes and relatively low elsewhere; $N_1$ points generated in step 1 of Algorithm \ref{alg:Self-normalized Sampling} can fail to have any point inside these large error volumes. A naive remedy to this is to generate more than $N_1$ points in step 1 so that it covers a broader range of points. On the other hand, the  Metropolis-Hastings algorithm can substantially sample points on these small volumes although the number of burn-ins needs to be very large to reach these volumes.

\begin{figure}
\centering
  \begin{minipage}[b]{0.45\textwidth}
    \caption{Distribution of Squared Residual Errors}\label{fig:3.1}
    \includegraphics[width=\textwidth]{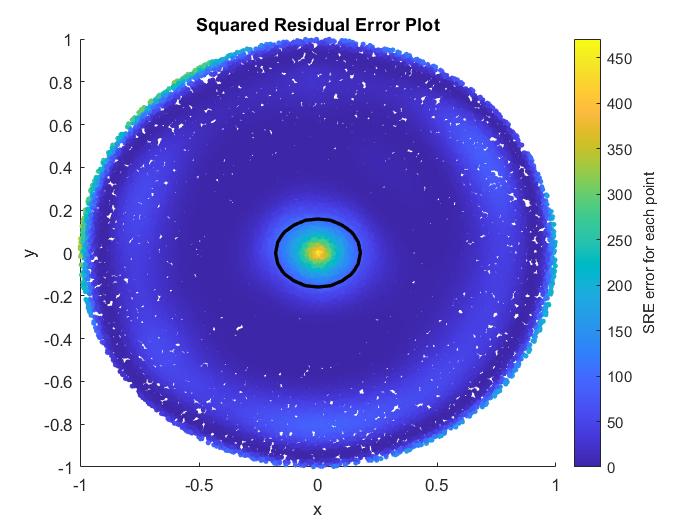}
  \end{minipage}
  \hfill
  \begin{minipage}[b]{0.45\textwidth}
    \caption{500 Points Generated by Uniform Annular}\label{fig:3.2}
    \includegraphics[width=0.89\textwidth]{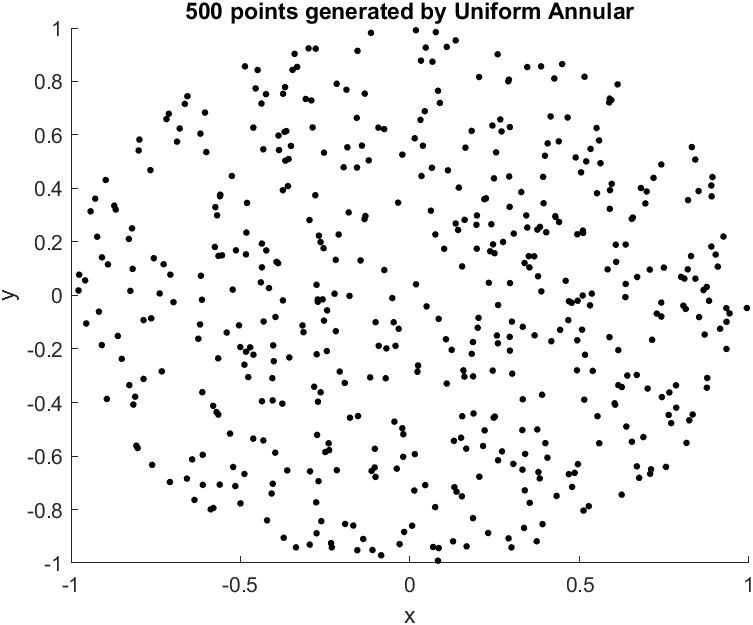}

  \end{minipage}
\end{figure}

\begin{figure}
\centering
  \begin{minipage}[!h]{0.45\textwidth}
    \caption{500 Points Generated by Metropolis}\label{fig:3.3}
    \includegraphics[width=\textwidth]{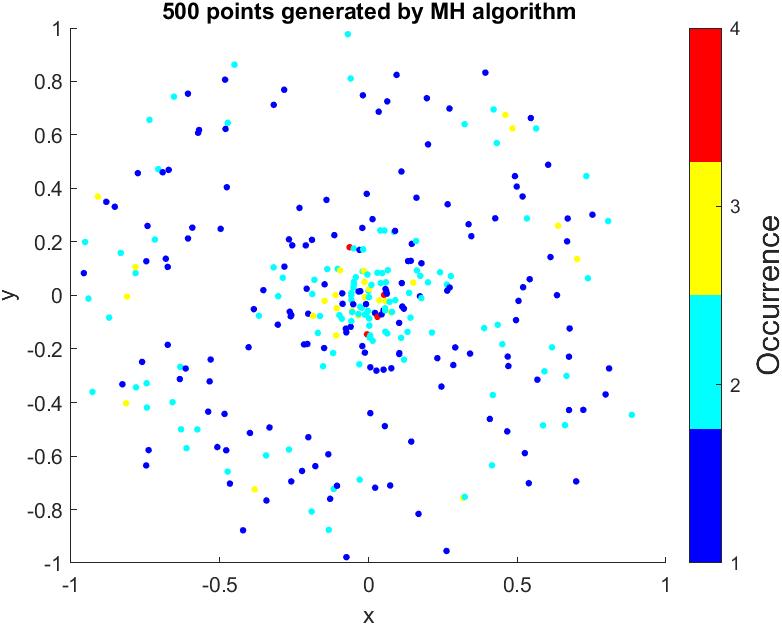}
  \end{minipage}
  \hfill
  \begin{minipage}[!h]{0.45\textwidth}
    \caption{500 Points Generated by Self-normalized}\label{fig:3.4}  
    \includegraphics[width=\textwidth]{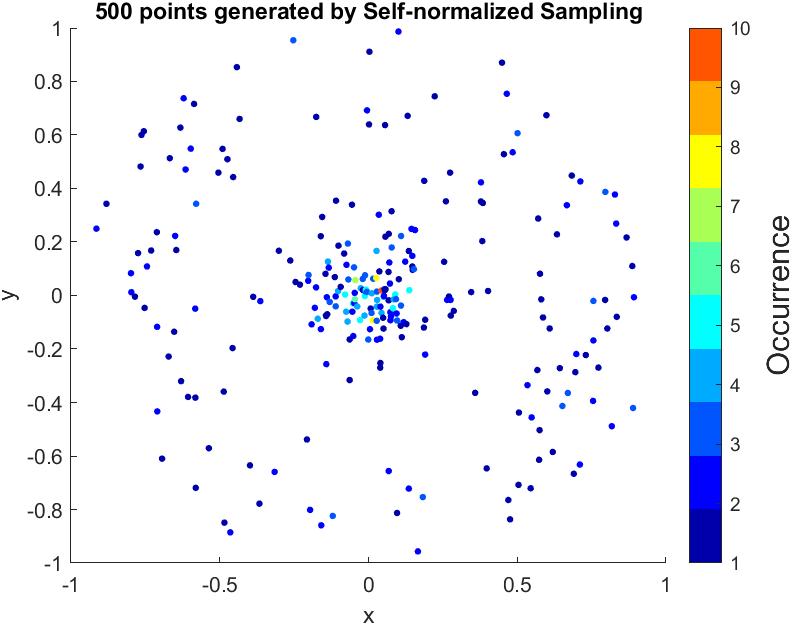}
  \end{minipage}
\end{figure}

A demonstration of points generated in 2D by Algorithms \ref{alg:Metropolis Hastings Sampling} and \ref{alg:Self-normalized Sampling} are presented. Figure \ref{fig:3.1}  is an example of a distribution of squared residual errors. Figures \ref{fig:3.3} and \ref{fig:3.4} are distributions of $500$ points generated by the Metropolis-Hastings algorithm and the self-normalized sampling, respectively, with $p = 1$ and $3500$ burn-ins for the Metropolis-Hastings algorithm;   the uniform distribution used in Algorithms \ref{alg:Metropolis Hastings Sampling} and  \ref{alg:Self-normalized Sampling} is replaced by uniform annular distribution  \cite{gu2021selectnet}. Figure \ref{fig:3.2} is the distribution of 500 points generated by a uniform annular distribution. The number of points inside of the ellipse $\frac{x^2}{0.18^2} + \frac{y^2}{0.16^2} = 1$ generated by each approach are shown in Table \ref{table:1}.

\begin{table}[h]
	\caption{The Number of Points Inside Ellipse}\label{table:1}
	\centering
	\begin{tabular}{ll}
		\toprule
		Approach & Number of Points \\
		\midrule
		Uniform Annular & 61   \\
		Metropolis Hastings & 210   \\
		Self-normalized Sampling & 255   \\
		\bottomrule
	\end{tabular}
\end{table}

\subsection{Adaptive Sampling: A Statistical Viewpoint}
Adaptive sampling is motivated by active learning that takes values of residual errors as a proxy of the model being wrong. It is natural to connect adaptive sampling with uncertainty sampling that takes uncertainty as a proxy of the model being wrong. Adaptive sampling can   be explained from a statistical viewpoint. The definition of expectation and Monte-Carlo estimate is in Appendix \ref{defs}.

\begin{definition}[Owen, 2013  \cite{mcbook}] 
Let $\boldsymbol{x} \in \Omega \subset \mathbb{R}^d$ be a random variable and $w(\boldsymbol{x})$ be any function defined on $\Omega$ with density $p(\boldsymbol{x})$. Let q(x) be another probability density on $\Omega$ such that $q(\boldsymbol{x}) = 0$ implies $w(\boldsymbol{x})p(\boldsymbol{x}) = 0$. $q(\boldsymbol{x})$ is also called biasing distribution  or important distribution. The importance sampling estimate is defined by
$$\hat{\mu}_{q}= \frac{|\Omega|}{n} \sum_{i=1}^{n} \frac{w\left(\boldsymbol{x}_i\right)p\left(\boldsymbol{x}_i\right) }{q\left(\boldsymbol{x}_i\right)} \text{, } \boldsymbol{x}_i \sim q(\boldsymbol{x}).$$
Denote by $Unif(\Omega)$ a uniform distribution on $\Omega$. Note that if $p(\boldsymbol{x}) = Unif(\Omega)$, then $|\Omega|$ and $p(\boldsymbol{x}) = \frac{1}{|\Omega|}$ cancel out each other. 
\end{definition}

\begin{thm}[Owen, 2013  \cite{mcbook}]\label{thm:3.1}
The importance sampling estimate $\hat{\mu}_{q}$  is an unbiased estimate of $\mu $.
\end{thm}

As one may have noticed, Theorem \ref{thm:3.1} can be applied to estimate expectations in Equation (\ref{eq:2.3}), where $|\mathcal{D} \phi(\boldsymbol{x} ;\boldsymbol{\theta})-f(\boldsymbol{x})|^{2} = w(\boldsymbol{x})$,  $q(\boldsymbol{x})$ is defined in (\ref{eq:3.2}), and $p(\boldsymbol{x}) = Unif(\Omega)$.  The  importance sampling estimate  of $\mu =  \mathbb{E}_{\boldsymbol{x} \in \Omega}\left[|\mathcal{D} \phi(\boldsymbol{x} ;\boldsymbol{\theta})-f(\boldsymbol{x})|^{2}\right] = \mathbb{E}_{\boldsymbol{x} \in \Omega} \left[w(\boldsymbol{x})\right]$ is thus
$ 
    \hat{\mu}_{q}=\frac{1}{n} \sum_{i=1}^{n} \frac{w\left(\boldsymbol{x}_i\right) }{q\left(\boldsymbol{x}_i\right)} \text{, } \boldsymbol{x}_i \sim q(\boldsymbol{x}) 
$.

One way to justify adaptive sampling is to adopt the notion of active learning. Adaptive sampling can also be understood as a better estimating approach. In  \eqref{eq:2.3}, the plain Monte-Carlo estimate  will face the issue of slow convergence if, for example, the distribution of squared residual errors is large at a small volume $\mathcal{V} \subset \Omega$ but is low, or nearly zero, outside of $\mathcal{V}$. The plain Monte-Carlo estimate could fail to have even one point inside the volume $\mathcal{V}$. Faster convergence can be achieved ``by sampling from a distribution that overweights the important region, hence the name importance sampling"  \cite{mcbook}. This paper is not concerned with finding an unbiased estimate to the expectation but finding the best parameter $\boldsymbol{\theta}$ that minimizes the expectation. Therefore, allocation points will be sampled from $q(\boldsymbol{x})$ without altering the formulation in the last step of   \eqref{eq:2.3} since only $w(\boldsymbol{x})$ depends on $\boldsymbol{\theta}$.

\section{Setup of Numerical Implementation}\label{sec:4}
\subsection{Network Architecture}
Note that adaptive sampling introduced in Section \ref{sec:3} is an active-learning-based sampling technique that is independent of the choice of network architectures. 
In the numerical implementation, the solution network $\phi(\boldsymbol{x}; \boldsymbol{\theta})$ is a feed-forward neural network in which each layer is fully connected. The network can be expressed as a composition of several activation functions; i.e.,   
\begin{equation}\label{eq:4.1}
\phi(\boldsymbol{x}; \boldsymbol{\theta}) = h_L\circ h_{L-1}\circ...\circ h_{1} \circ h_{0}(\boldsymbol{x}),
\end{equation}
 
$$h_l(\boldsymbol{x}^l) = \sigma(\boldsymbol{W}^l \boldsymbol{x}^l + b^l), \text{ for } l = 0, 1, 2, ..., L-1,\qquad h_L(\boldsymbol{x}^L) = \boldsymbol{W}^L (\boldsymbol{x}^L + b^L),$$
where $L$ is  the depth of the network, $\sigma$ is a nonlinear activation function,   
 $\boldsymbol{W}^{0} \in \mathbb{R}^{m \times d}$,   $\boldsymbol{W}^{L} \in \mathbb{R}^{m \times 1}$,  
$\boldsymbol{W}^{l} \in \mathbb{R}^{m \times m}$, $l=1, \cdots, L-1$, $\boldsymbol{b}^{L} \in \mathbb{R}$, $\boldsymbol{b}^{l} \in \mathbb{R}^{m \times 1}$, $l=0, \cdots, L-1$, $m$ is the number of nodes or neurons in a layer and called the width of the network,  $d$ is the dimension of the problem and  $\boldsymbol{\theta} =  \{\boldsymbol{\boldsymbol{W}}^{l}$, $\boldsymbol{b}^{l}\}_{l=0}^{L}$. Network setting parameters used in the numerical implementation are listed in Table \ref{table:2}.
\begin{table}[!h]
	\caption{The Setting of the Solution Network}\label{table:2}
	\centering
	\begin{tabular}{lll}
		\toprule
		Symbol & Meaning & Value\\
		\midrule
		$L$ & Depth of the Network & 3    \\
		$m$ & Width of the Network & 100   \\
		$\sigma$ & Activation Function& $ReLU^3$ = $max{(x^3, 0)}$\\
        $d$ & Dimension of Spatial Coordinates & to be specified\\
        $\tau$ & Learning Rate for Solution Network $\phi$ & Specified in \eqref{eq:lr}\\
		\bottomrule
	\end{tabular}
\end{table}

The network is trained to solve the minimization problem defined in (\ref{eq:2.3}) via Adam \cite{adam} (a variant of stochastic gradient descent), where the empirical loss is defined by
\begin{equation}\label{eq:4.3}
J(\boldsymbol{\theta}) = \frac{1}{N_{1}} \sum_{i=1}^{N_{1}}\left|\mathcal{D} \phi\left(\boldsymbol{x}_{i} ; \boldsymbol{\theta}\right)-f\left(\boldsymbol{x}_{i}\right)\right|^{2}+ \frac{\lambda}{N_{2}} \sum_{j=1}^{N_{2}}\left|\mathcal{B} \phi\left(\boldsymbol{x}_{j} ; \boldsymbol{\theta}\right)-g\left(\boldsymbol{x}_{j}\right)\right|^{2},
\end{equation}
and $\boldsymbol{\theta}$ is updated by 
\begin{equation}
    \boldsymbol{\theta} \leftarrow \boldsymbol{\theta} - \tau \nabla J(\boldsymbol{\theta}),
\end{equation}
where $\tau$ is the learning rate of the solution network  defined by 
\begin{equation}\label{eq:lr}
\begin{aligned}
    &\tau^{(k)} = 10^{-3-(3j/1000)}, \text{ for } \ceil[\Big]{\frac{0.999n}{1000}j} \leq k < \ceil[\Big]{\frac{0.999n}{1000}(j+1)}, \text{ } j = 0, 1, 2, ..., 999,\\
    &\tau^{(k)} = 10^{-6}, \text{ for } \ceil[\Big]{0.999n} \leq k < n,
\end{aligned}
\end{equation}
where $\tau^{(k)}$ denotes the learning rate at the $k$-th epoch, and $n$ is the number of total epochs. Index starts from 0, i.e., the learning rate for the first epoch is $\tau^{(0)}$. 

\subsection{Derivatives of Networks}
In order to calculate residual errors, $\mathcal{D}\phi(\boldsymbol{x}; \boldsymbol{\theta})$ needs to be evaluated. It is well-known that Autograd can perform such differentiation. However, in our experiment, we observe that Autograd does not output precise partial derivatives in high dimensions. Therefore, in our implementation, the numerical differentiation (see Appendix \ref{num_diff}) with step size $h = 0.0001$ will be used to estimate partial derivatives except in Section \ref{sec:5.4},  where Autograd is used since it outputs correct partial derivatives in 2D. For example, let $\phi(\boldsymbol{x}; \boldsymbol{\theta})$ be the solution network with $\boldsymbol{x}\in \mathbb{R}^d$ and parameters $\boldsymbol{\theta}$.  The numerical differentiation  estimate of $\phi(\boldsymbol{x}; \boldsymbol{\theta})$, denoted by $\nabla \phi(\boldsymbol{x}; \boldsymbol{\theta})$, is defined by 
$$\nabla \phi(\boldsymbol{x}; \boldsymbol{\theta}) \approx \frac{1}{h} \sum_{i = 1}^d \big[ \phi(\textbf{a} + h\boldsymbol{e}_i; \boldsymbol{\theta}) - \phi(\textbf{a};\boldsymbol{\theta}) \big],$$
where $\textbf{a} \in \mathbb{R}^d$,   $\boldsymbol{e}_i$ is a vector of all $0$'s except a $1$ in the i-th entry, and $h\in \mathbb{R}$ is the step size.
Note that, as found in \cite{gu2021selectnet}, with a step size of $10^{-4}$, the truncation error is up to $O(10^{-6})$, and it can be ignored in practice since the final error is at least $O(10^{-4})$. We will  resolve the issue that Autograd does not output precise partial derivatives in our future work. As mentioned above, numerical differentiation only causes negligible errors in this work.

\subsection{Training Solution Networks with Active Sampling}
\begin{algorithm}[H]
\caption{Training Solution Networks}\label{alg:4.1}
\SetAlgoLined
\setcounter{AlgoLine}{0}
\SetKwInput{KwResult}{Result}
\SetKwRepeat{Do}{do}{while}

\KwResult{parameters $\boldsymbol{\theta}^\ast$}
\SetKwInput{KwRequire}{Require}
\KwRequire{PDE (\ref{eq:2.1})}
 Set $n = $ total iterations/epoch, $N_1$ and $N_2$ for size of sample points in $\Omega$ and on $\partial \Omega$ respectively 
 \;
 Initialize $\phi(\boldsymbol{x}; \boldsymbol{\theta})$;
 \While{$k < n$}{
    Generate sampling training points by Algorithm \ref{alg:Metropolis Hastings Sampling} or \ref{alg:Self-normalized Sampling},
    $\left\{\boldsymbol{x}_{i}^{1}\right\}_{i=1}^{N_{1}} \subset \Omega \text { and }\left\{\boldsymbol{x}_{j}^{2}\right\}_{j=1}^{N_{2}} \subset \partial \Omega$ 
    \;
    Loss: $J(\boldsymbol{\theta}) = \frac{1}{N_{1}} \sum_{i=1}^{N_{1}}\left|\mathcal{D} \phi\left(\boldsymbol{x}^1_{i} ; \boldsymbol{\theta}\right)-f\left(\boldsymbol{x}^1_{i}\right)\right|^{2}+ \frac{\lambda}{N_{2}} \sum_{j=1}^{N_{2}}\left|\mathcal{B} \phi\left(\boldsymbol{x}^2_{j} ; \boldsymbol{\theta}\right)-g\left(\boldsymbol{x}^2_{j}\right)\right|^{2}$
    \;
    Update: $\boldsymbol{\theta} := \boldsymbol{\theta} - \tau \nabla J(\boldsymbol{\theta})$
    \;
    }
 \Return $\boldsymbol{\theta}^\ast := \boldsymbol{\theta}$ \;
\end{algorithm}

\subsection{Uniform Annular Distribution}
The uniform annular approach will be employed to sample points in $\Omega$ in the first step of Algorithm \ref{alg:Metropolis Hastings Sampling} and \ref{alg:Self-normalized Sampling}. The uniform annular approach divides the interior of the domain  $\Omega$  into $N_{a}$ annuli $\left\{k / N_{a}<|\boldsymbol{x}|<\right.$ $\left.(k+1) N_{a}\right\}_{k=0}^{N_{a}-1}$ and generates $N_{1} / N_{a}$ samples uniformly in each annulus, where   $N_1$ is the number of training points in $\Omega$ in each epoch, and $N_a$ is the number of annuli; in practice $N_a$ should be divided by $N_1$, i.e., $N_a \vert N_1$. This distribution is also utilized in \cite{gu2021selectnet}. As a side note, this uniform annular distribution approach can also be understood as a biased distribution that covers a broader range of points in $\Omega$. Unlike the adaptive sampling distribution in  (\ref{eq:3.2}), the uniform annular distribution is not necessarily a better distribution for the model to learn for each PDE example. It is worth pointing out that the domain $\Omega$ considered in this paper is always a
high-dimensional unit ball or unit cube. Therefore, this uniform annular distribution can be employed. In the general case, it is difficult to apply the uniform annular distribution.
\subsection{Accuracy Assessment}
To measure the accuracy of the model, 10000 testing points $\left\{\boldsymbol{x}_{i}^{t}\right\}_{i=1}^{10000} \subset \Omega$, which are different from training points, will be sampled from the uniform annular distribution to approximate the overall relative $\ell_2$ error and the relative $\ell_1$ maximum modulus error at these points. Note that in the implementation, the random seed is fixed which implies that for the same PDE example, all models will be tested on the same set of uniform annular points to measure accuracy. The overall relative $\ell_2$ error and the relative $\ell_1$ maximum modulus error are, respectively, defined as follows 
\begin{equation*}
    e^{overall}_{\ell^{2}}(\boldsymbol{\theta}):=\frac{\left(\sum_{i = 1}^{10000}\left|\phi\left(\boldsymbol{x}_{i}^t ; \boldsymbol{\theta}\right)-u\left(\boldsymbol{x}_{i}^t\right)\right|^{2}\right)^{\frac{1}{2}}}{\left(\sum_{i = 1}^{10000}\left|u\left(\boldsymbol{x}_{i}^t\right)\right|^{2}\right)^{\frac{1}{2}}},
    \quad \quad
    e^{max}_{modulus}(\boldsymbol{\theta}):=\frac{max\big(\left|\phi\left(\boldsymbol{x}_{i}^t ; \boldsymbol{\theta}\right)-u\left(\boldsymbol{x}_{i}^t\right)\right|\big)_{i = 1}^{10000}}{max\big(\left|u\left(\boldsymbol{x}_{i}^t\right)\right|\big)_{i = 1}^{10000}}.
\end{equation*}

\section{Numerical Experiment}\label{sec:5}
 In this section, our adaptive sampling is tested on three types of PDE examples: elliptic, parabolic, and hyperbolic equations. Testing includes various dimensions. In this section, we will always use Algorithm \ref{alg:Self-normalized Sampling} (Self-normalized Sampling) to sample training points. In our experiments, Algorithm \ref{alg:Metropolis Hastings Sampling} (Metropolis Hastings Sampling) and Algorithm \ref{alg:Self-normalized Sampling} (Self-normalized Sampling) give similar results. It is worth mentioning that the main focus of this work is to ``choose'' training examples based on the residual error distribution, and it is independent of the algorithm that is used to simulate the residual error distribution. Therefore, one may use other algorithms, such as Gibbs sampling, to simulate the distribution as well. 
Adaptive sampling will be compared with RAR in Examples \ref{sec:elliptic}, \ref{sec:parabolic}, and \ref{sec: hyperb}. In the comparison experiments, we only replace AS with RAR (in line $3$ of Algorithm \ref{alg:4.1}, AS is replaced by RAR), and all other settings are kept the same. The RAR algorithm used in this work (see Algorithm \ref{alg:6.1}) is slightly different from the  algorithm  proposed in the paper \cite{deep_XDE} in order to make the comparison in which the number of training points will still be fixed. $m$ in the algorithm proposed in the paper \cite{deep_XDE} is chosen to be $2000$ for testing. 
Finally, our adaptive sampling will be applied to some current frameworks to test if adaptive sampling is compatible with these frameworks. The parameter setting for Sections 4.1 and 4.2 is listed in Table \ref{table:3}. 

\begin{table}[h!]
	\caption{Parameter Setting}\label{table:3}
	\centering
	\begin{tabular}{lll}
		\toprule
		Symbol & Meaning & Value\\
		\midrule
		$n$ & Number of Epochs & 20000    \\
		$N_1$ & Number of Training Points in $\Omega$ in Each Epoch & 12000  \\
		$N_2$ for Example \eqref{eq4.1} &  Number of Training Points on $\partial \Omega$ in Each Epoch& 12000\\
		$N_2$ for Example \eqref{eq4.2} &  Number of Training Points on $\partial \Omega$ in Each Epoch& 12000 + $\frac{12000}{d}$\\
		$N_2$ for Example \eqref{eq:hyperbo} &  Number of Training Points on $\partial \Omega$ in Each Epoch & 12000 + $\frac{12000}{d}$\\

        $\lambda$ & Boundary Loss Weighting Term in \eqref{eq:hyperbo} & 10\\
		\bottomrule
	\end{tabular}
\end{table}

\subsection{Elliptic Equation}\label{sec:elliptic}
\noindent We consider a nonlinear elliptic equation:
\begin{equation}\label{eq4.1}
\begin{aligned}
-\nabla \cdot \left((1+\frac{1}{2}|\boldsymbol{x}|^{2}) \nabla u\right)+(\nabla u)^{2} &=f(\boldsymbol{x}),  \quad \text { in } \Omega:=\{\boldsymbol{x}:|\boldsymbol{x}|<1\}, \\
u &=g(\boldsymbol{x}), \quad  \text { on } \partial \Omega,
\end{aligned}
\end{equation}
where $g(\boldsymbol{x}) = 0$ and $f$ is given appropriately so that the exact solution is given by $
u(\boldsymbol{x})= sin(\frac{\pi}{2}(1-|\boldsymbol{x}|)^{2.5})$.
 
Adaptive sampling is tested on this PDE example in various high dimensions: 10, 20, and 100. Adaptive sampling will be compared with the  basic least squares method discussed in Section \ref{sec:2.1}. Since the computational costs of these two models are different, for this particular example, error decay versus time in seconds is also displayed to demonstrate the efficiency of adaptive sampling. Running time is obtained by training on Quadro RTX 8000 Graphics Card. Only training time has been taken into account; time for accuracy assessment and saving data is not included. 
The overall relative error seems very large in 100 dimensions because the true solution vanishes in most of the volume when the dimension is high.
\begin{table}
	\caption{Example \ref{sec:elliptic} Result; Error Reduction Formula: $1-\frac{error(AS)}{error(Basic)}$}\label{table:4}
	\centering
	\begin{tabular}{llllll}
		\toprule
		Dimension& &  AS & Basic & RAR & Error Reduced by AS\\
		\midrule
		 10  & \begin{tabular}{@{}l@{}@{}}
                   $\ell_2$ error\\
                   max $\ell_1$ error\\
                   time in Sec\\
                   \\
                 \end{tabular}& \begin{tabular}{@{}l@{}@{}}
                   8.784735e-03\\
                   3.681198e-02\\
                   8179.919345\\

                   \\
                 \end{tabular}  &\begin{tabular}{@{}l@{}@{}}
                   2.526952e-02\\
                   1.336612e-01\\
                   6155.072389\\
                   \\
                 \end{tabular} &\begin{tabular}{@{}l@{}@{}}
                   2.024158e-02\\
                   9.741917e-02\\
                   7423.530141\\
                   \\
                 \end{tabular} &\begin{tabular}{@{}l@{}@{}}
                   65.24$\%$\\
                   72.46$\%$\\
                   \\
                   \\
                 \end{tabular}\\

		 20 & \begin{tabular}{@{}l@{}@{}}
                   $\ell_2$ error\\
                   max $\ell_1$ error\\
                   time in Sec\\
                   \\
                 \end{tabular}& \begin{tabular}{@{}l@{}@{}}
                   3.102093e-02\\
                   1.145958e-01\\
                   12489.657977\\
                   
                   \\
                 \end{tabular}  &\begin{tabular}{@{}l@{}@{}}
                   7.198276e-02\\
                   2.719710e-01\\
                   7985.964646\\
                   \\
                 \end{tabular} &\begin{tabular}{@{}l@{}@{}}
                   6.345444e-02\\
                   3.157725e-01\\
                   9022.221589\\
                   \\
                 \end{tabular}&\begin{tabular}{@{}l@{}@{}}
                   56.90$\%$\\
                   57.86$\%$\\
                   \\
                   \\
                 \end{tabular}\\

		 100 & \begin{tabular}{@{}l@{}@{}}
                   $\ell_2$ error\\
                   max $\ell_1$ error\\
                   time in Sec\\
                   \\
                 \end{tabular}& \begin{tabular}{@{}l@{}@{}}
                   1.205655e-01\\
                   3.280543e-01\\
                   43374.188778\\
                   
                   \\
                 \end{tabular}  &\begin{tabular}{@{}l@{}@{}}
                   3.964277e-01\\
                   1.246953e+00\\
                   32276.845542\\
                   \\
                 \end{tabular}    &\begin{tabular}{@{}l@{}@{}}
                   3.423974e-01\\
                   1.016873e+00\\
                   37841.390514\\
                   \\
                 \end{tabular}&\begin{tabular}{@{}l@{}@{}}
                   69.59$\%$\\
                   73.69$\%$\\
                   \\
                   \\
                 \end{tabular}\\
		\bottomrule
	\end{tabular}
\end{table}

\begin{figure} \caption{Example \ref{sec:elliptic} numerical results in 10 dimensions}\label{fig:5.1}
\minipage{0.45\textwidth}
  \includegraphics[width=\linewidth]{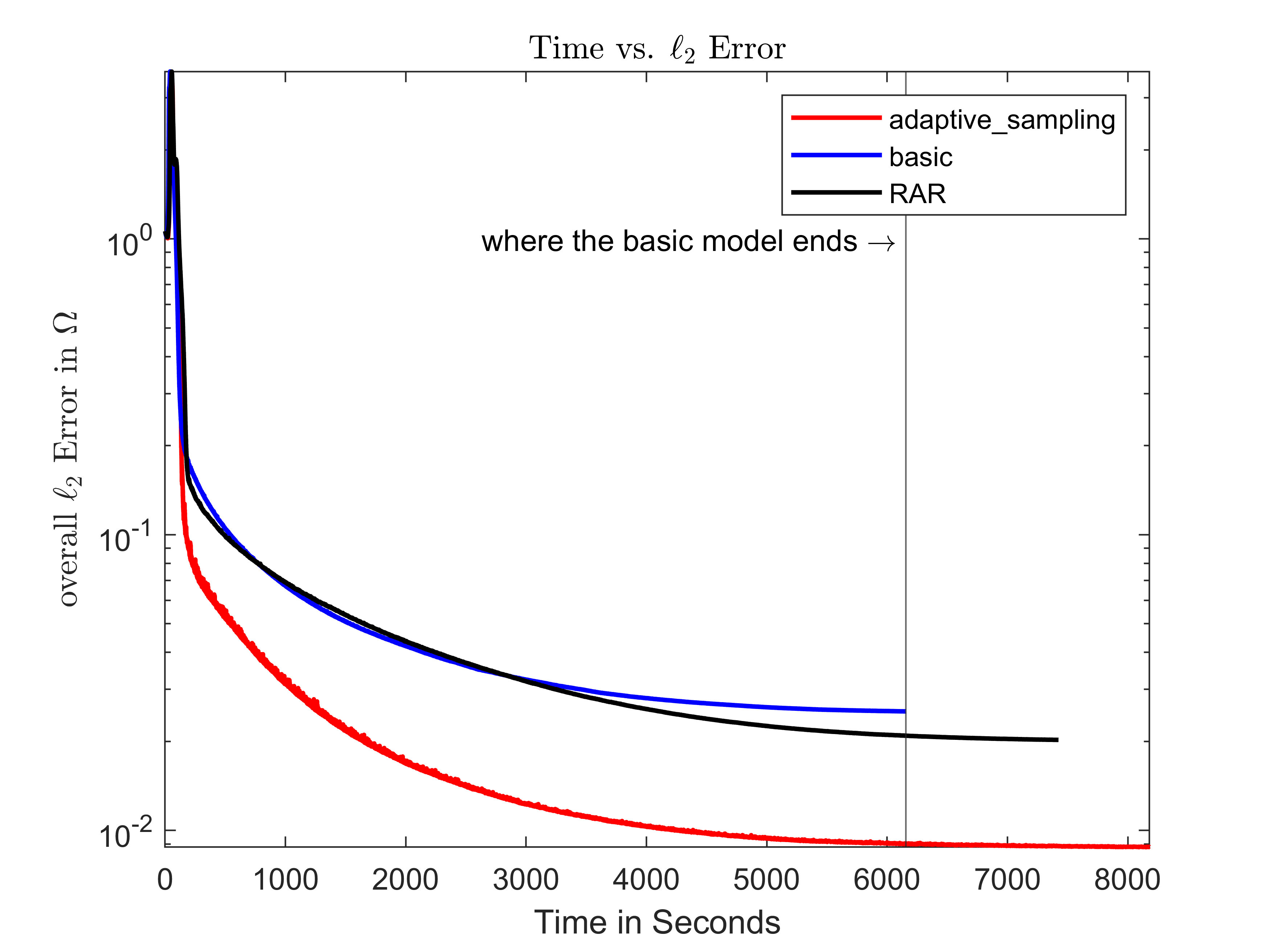}
  \caption*{(a) Time vs. Overall $\ell_2$ Error in $\Omega$}
\endminipage\hfill
\minipage{0.45\textwidth}
  \includegraphics[width=\linewidth]{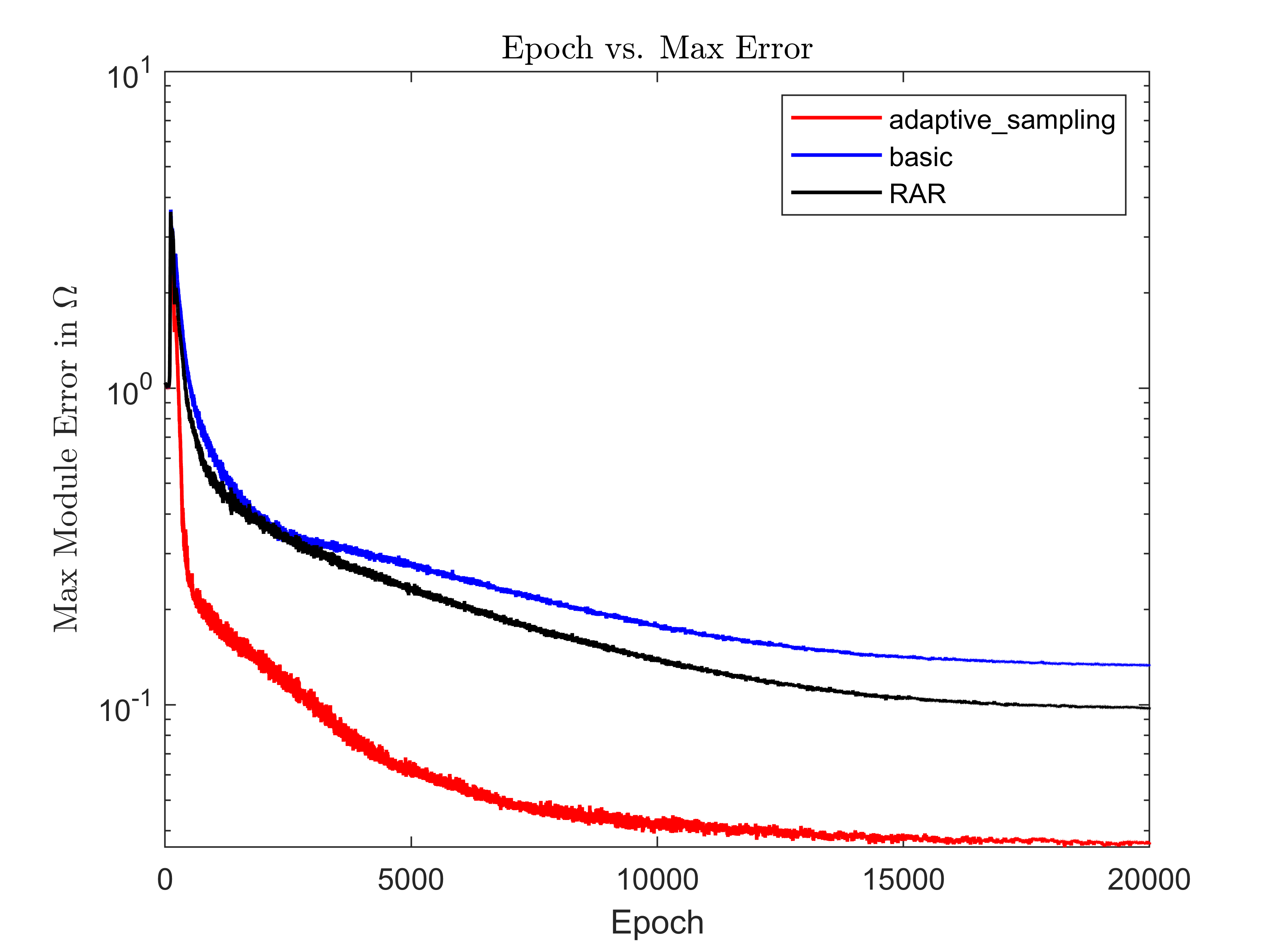}
  \caption*{(b) Epoch vs. Max Modulus Error in $\Omega$}
\endminipage\hfill
\end{figure}

\begin{figure} \caption{Example \ref{sec:elliptic} 10 dimensions $(x_1, x_2, 0, 0,...,0)$-surface of network solutions and the true solution}\label{fig:5.2}
\minipage{0.33\textwidth}
  \includegraphics[width=\linewidth]{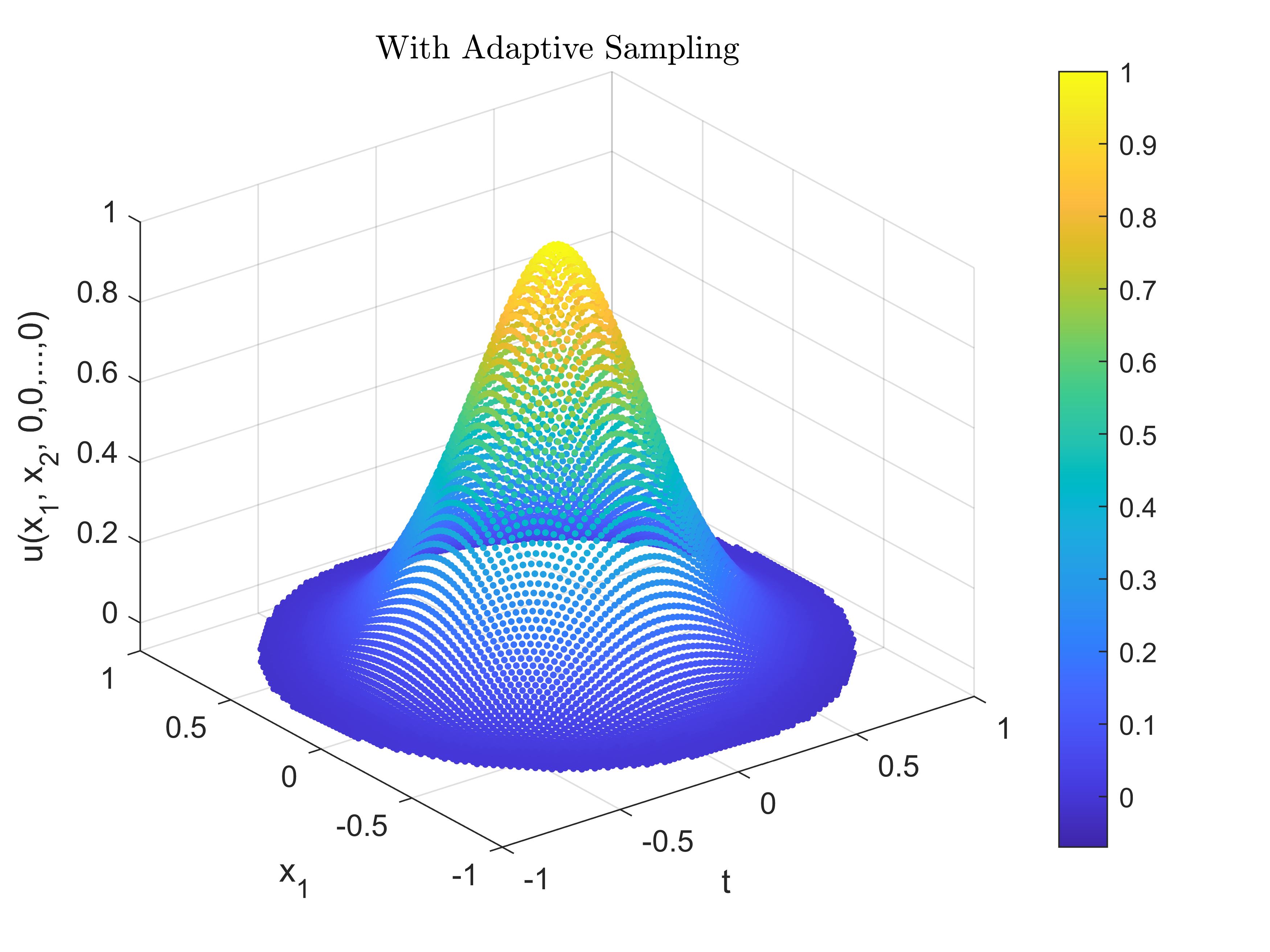}
  \caption*{(a) With Adaptive Sampling}
\endminipage\hfill
\minipage{0.33\textwidth}
  \includegraphics[width=\linewidth]{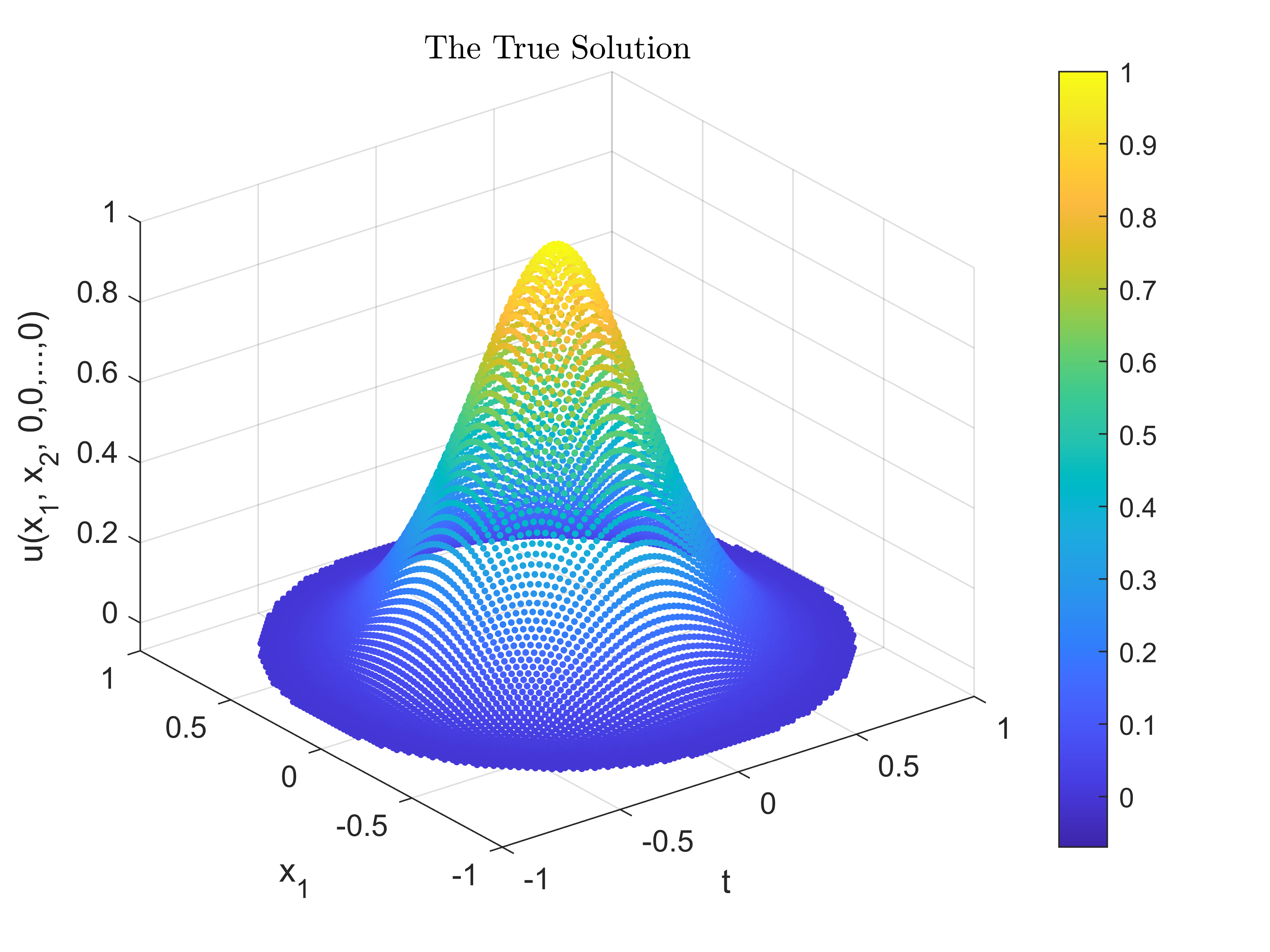}
  \caption*{(b) The True Solution}
\endminipage\hfill
\minipage{0.33\textwidth}
  \includegraphics[width=\linewidth]{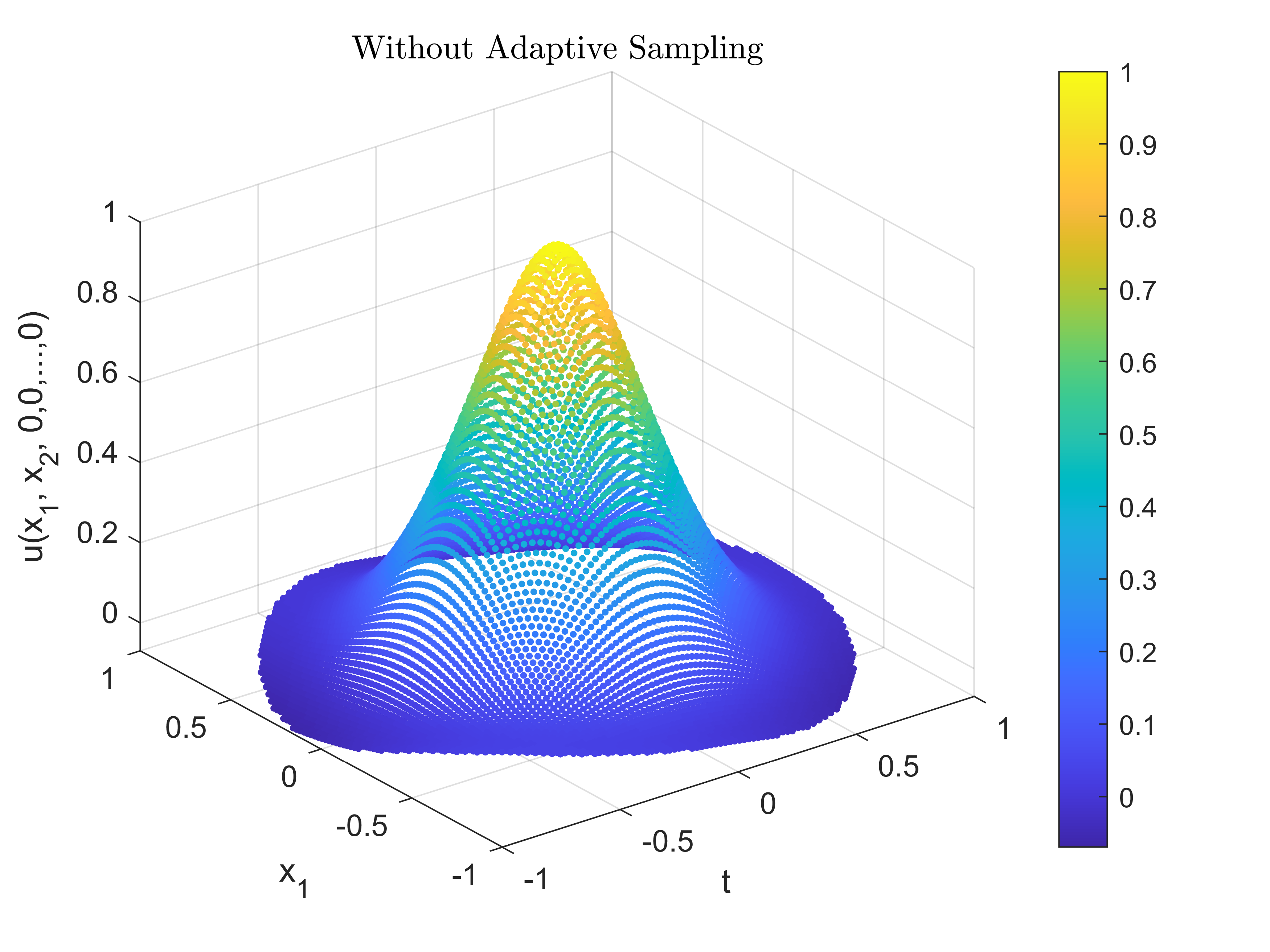}
  \caption*{(c)  Without  Adaptive Sampling}
\endminipage\hfill
\end{figure}

\begin{figure} \caption{Example \ref{sec:elliptic} 10 dimensions $(x_1, x_2, 0, 0,...,0)$-surface absolute difference $|u - \phi|$}\label{fig:5.3}
\minipage{0.45\textwidth}
  \includegraphics[width=\linewidth]{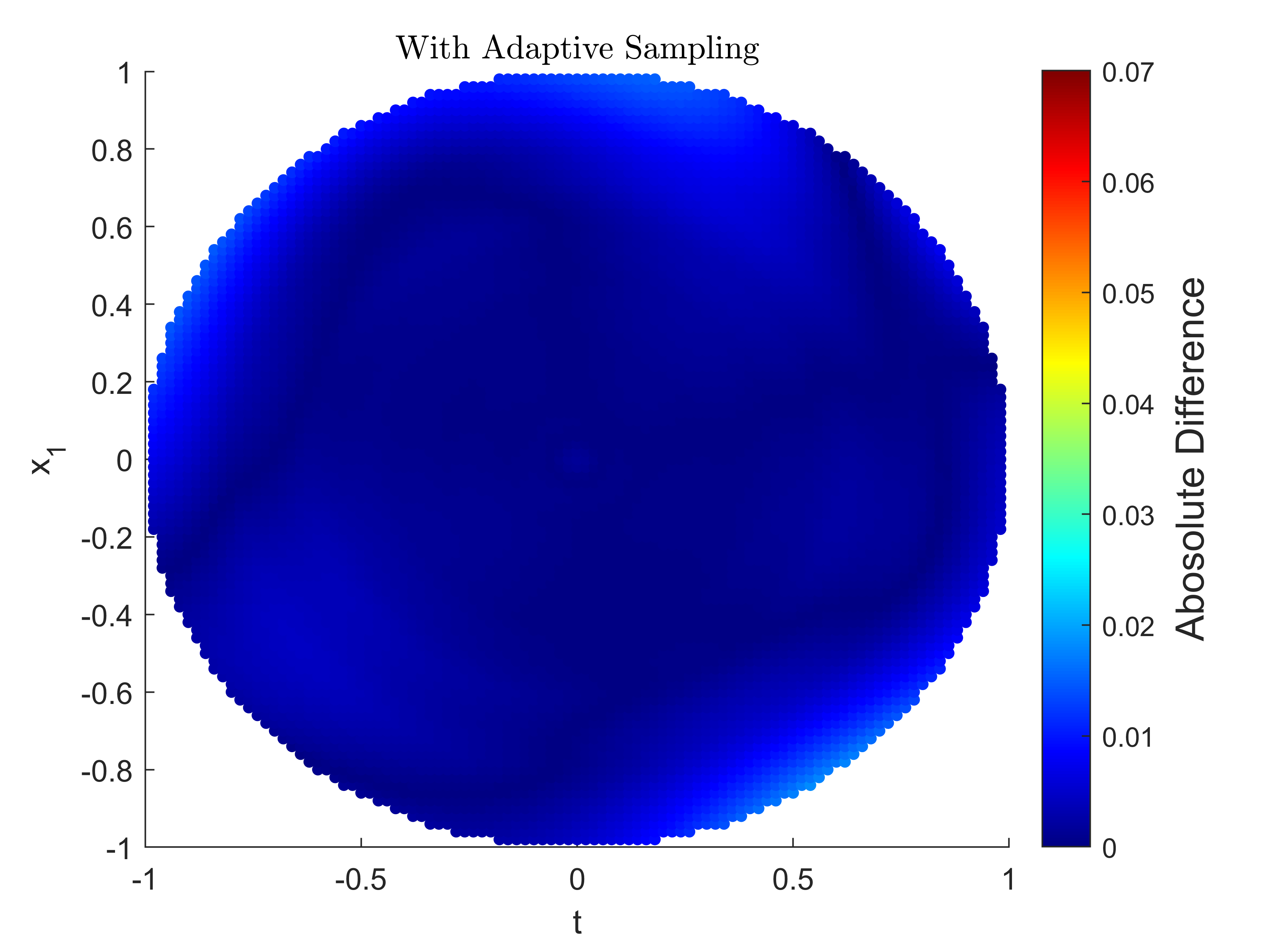}
  \caption*{(a) With Adaptive Sampling}
\endminipage\hfill
\minipage{0.45\textwidth}
  \includegraphics[width=\linewidth]{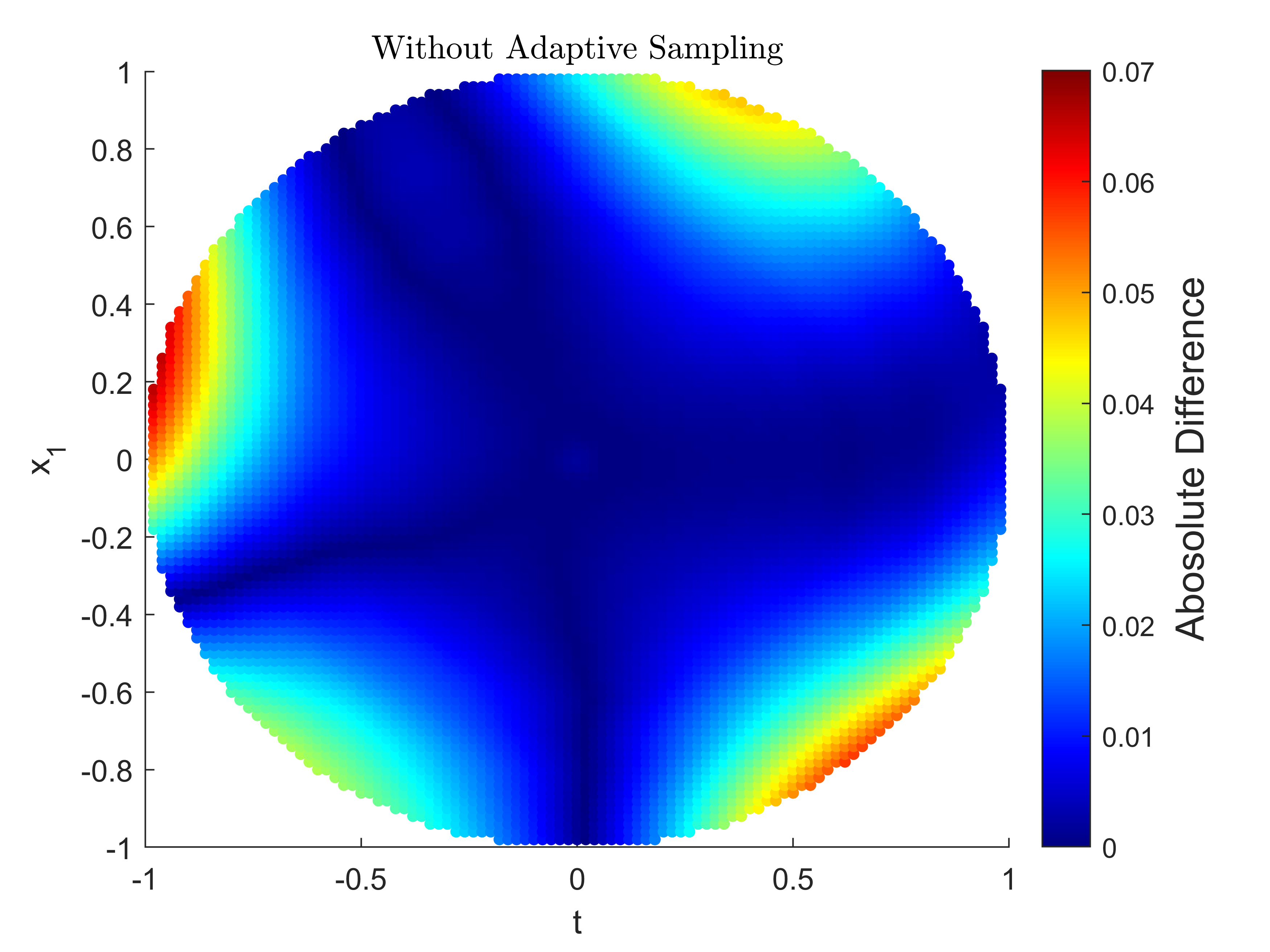}
  \caption*{(b) Without Adaptive Sampling}
\endminipage\hfill
\end{figure}

\begin{figure} \caption{Example \ref{sec:elliptic}  numerical results in 20 dimensions}\label{fig:5.4}
\minipage{0.45\textwidth}
  \includegraphics[width=\linewidth]{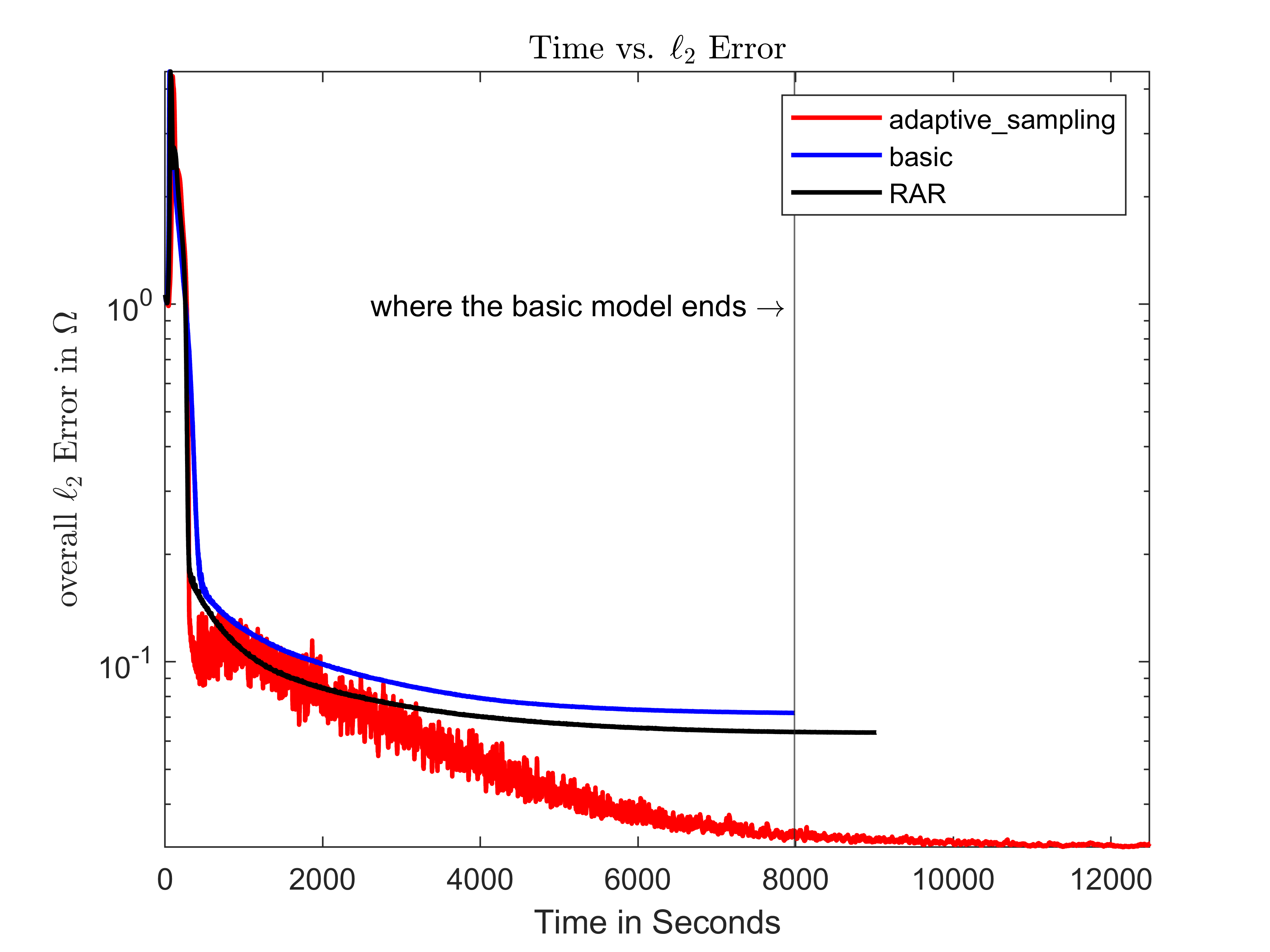}
  \caption*{(a) Time vs. Overall $\ell_2$ Error in $\Omega$}
\endminipage\hfill
\minipage{0.45\textwidth}
  \includegraphics[width=\linewidth]{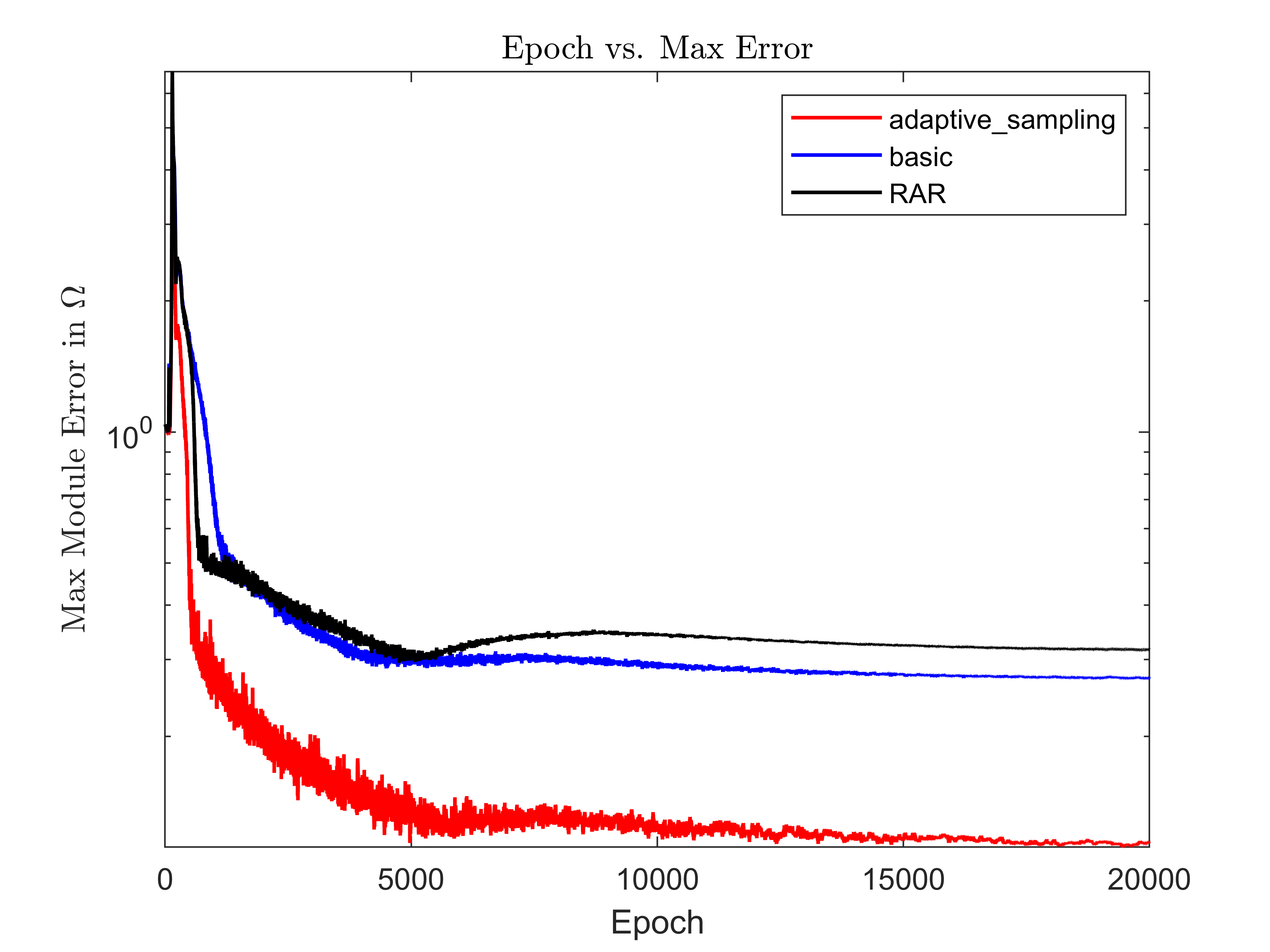}
  \caption*{(b) Epoch vs. Max Modulus Error in $\Omega$}
\endminipage\hfill
\end{figure}

\begin{figure} \caption{Example \ref{sec:elliptic} numerical results in 100 dimensions}\label{fig:5.5}
\minipage{0.45\textwidth}
  \includegraphics[width=\linewidth]{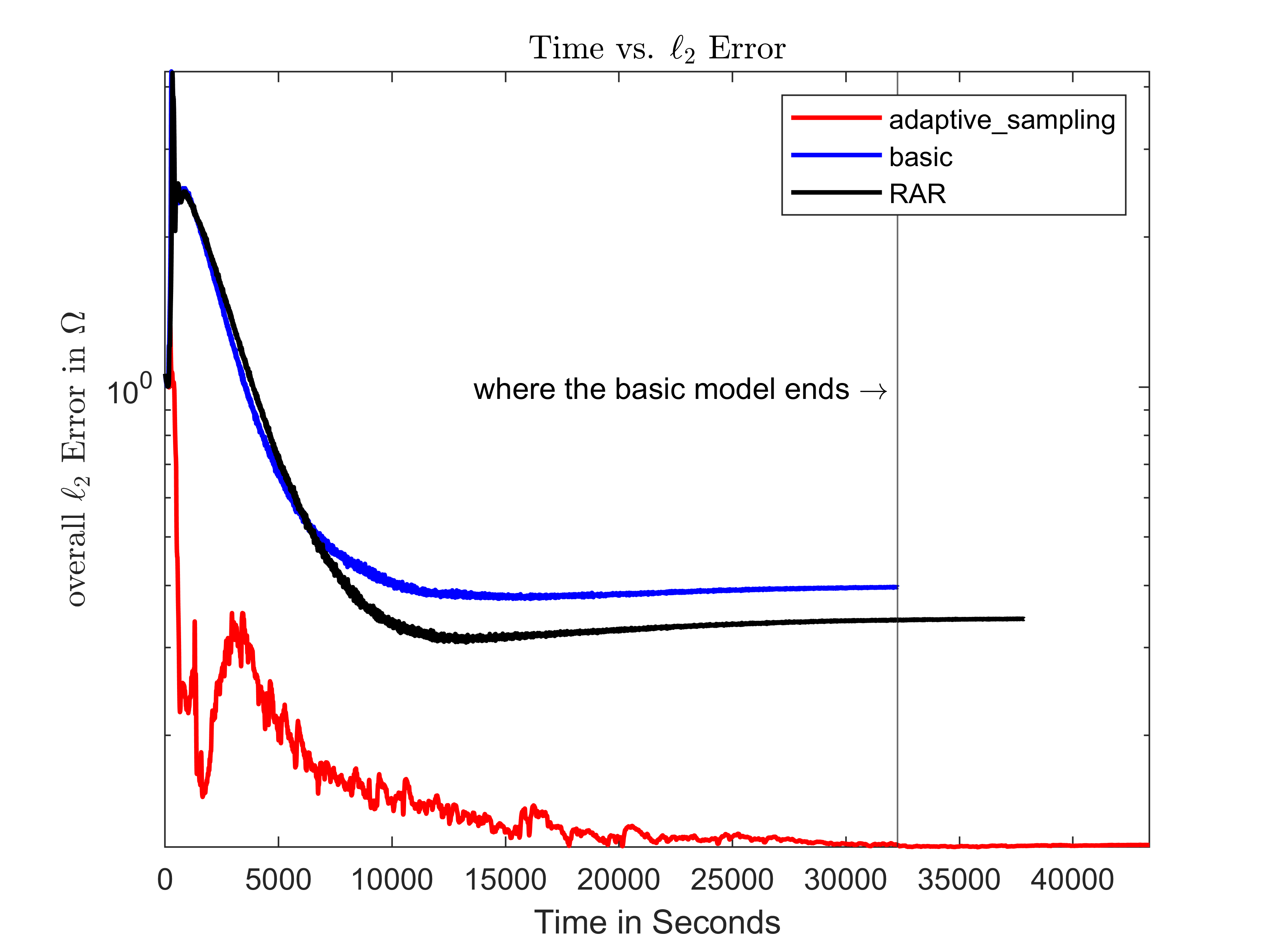}
  \caption*{(a) Time vs. Overall $\ell_2$ Error in $\Omega$}
\endminipage\hfill
\minipage{0.45\textwidth}
  \includegraphics[width=\linewidth]{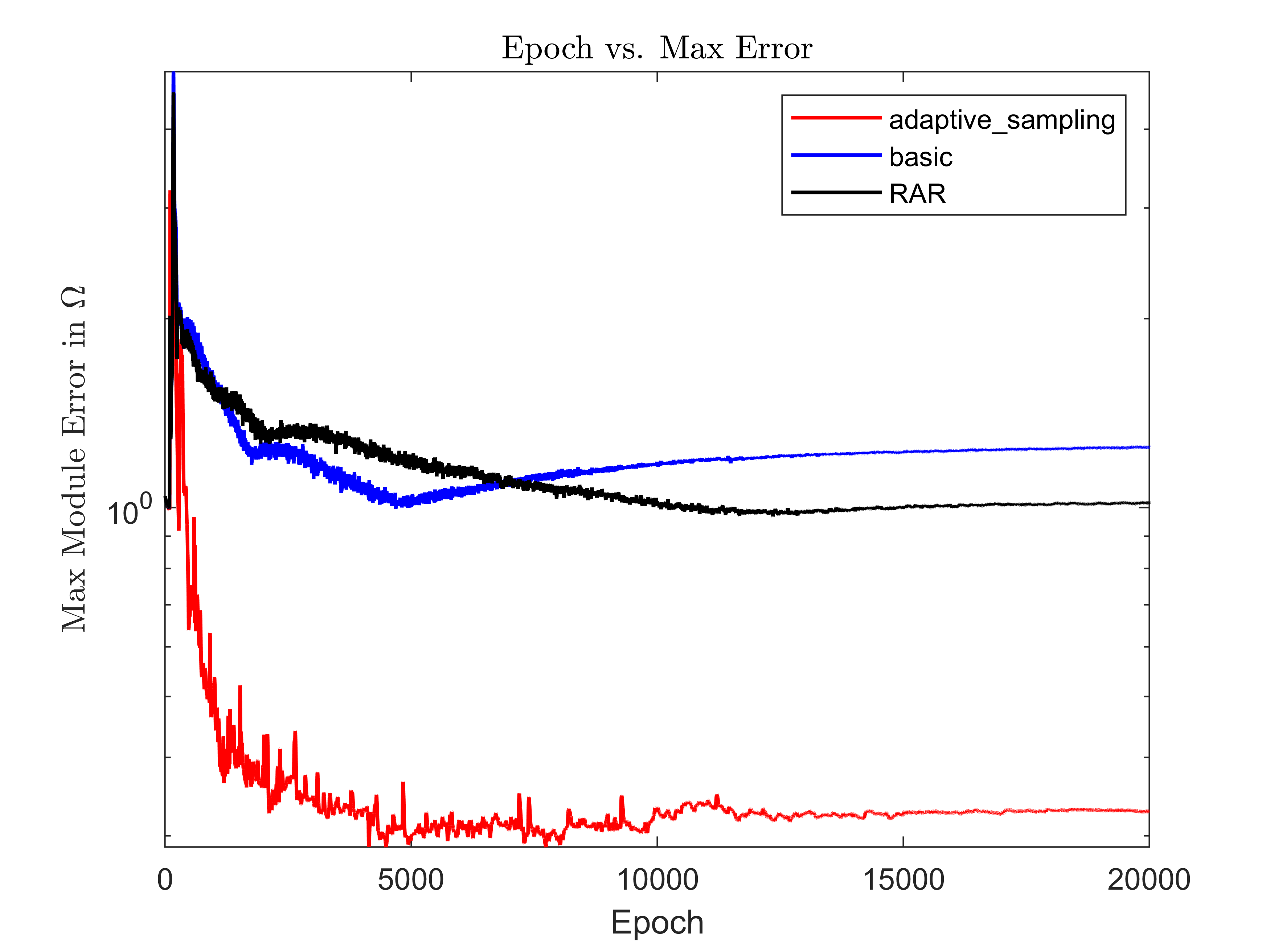}
  \caption*{(b) Epoch vs. Max Modulus Error in $\Omega$}
\endminipage\hfill
\end{figure}

Table \ref{table:4} shows the overall relative $l_2$ error and maximum relative $l_1$ modulus error by adaptive sampling and the  basic residual model  at the end of 20000 epochs. Time versus $\ell_2$ error decay in 10 dimensions, 20 dimensions, and 100 dimensions are demonstrated in Figures \ref{fig:5.1}, \ref{fig:5.4}, \ref{fig:5.5}, respectively. 
Figure \ref{fig:5.2} shows the $(x_1, x_2)$-surface of the ground-truth solution and network solutions with and without adaptive sampling, where axes are $x_1$, $x_2$, and $u(x_1, x_2,0, 0, ..., 0)$. Figure \ref{fig:5.3} shows the heat-map of the absolute difference in $(x_1, x_2)$-surface, where the color bar represents the absolute difference between the ground-truth solution and network solution. These results clearly show that adaptive sampling helps to reduce the error, especially relative $l_1$ maximum modulus error. More precisely, adaptive sampling helps to significantly reduce the variance in the distribution of error over the domain so that, compared to the results without adaptive sampling, it is less likely to have volumes/areas where the error is relatively much larger. 

\subsection{Parabolic Equation}\label{sec:parabolic}
\noindent We consider the following parabolic equation:
\begin{equation}\label{eq4.2}
\begin{aligned}
\partial_{t} u(\boldsymbol{x}, t)-\nabla_{\boldsymbol{x}} \cdot\left((1+\frac{1}{2}|\boldsymbol{x}|) \nabla_{\boldsymbol{x}} u(\boldsymbol{x}, t)\right) =f(\boldsymbol{x}, t)&, \quad  \text {in } \Omega:=\omega \times \mathbb{T}, \\
u(\boldsymbol{x}, t) =g(\boldsymbol{x}, t)&, \quad  \text {on } \partial \Omega = \partial \omega \times \mathbb{T},\\
u(\boldsymbol{x}, 0) =h(\boldsymbol{x})&, \quad \text {in } \omega,
\end{aligned}
\end{equation}
  where $\omega:=\{\boldsymbol{x}:|\boldsymbol{x}|<1\}$, $\mathbb{T} = (0, 1)$
, and 
$$
g(\boldsymbol{x}) = e^{|\boldsymbol{x}|\sqrt{1-t}},
$$
$$
h(\boldsymbol{x}) = exp(|\boldsymbol{x}|),
$$
$f$ is given approriately so that the exact solution is given by 
$u(\boldsymbol{x}, t) = e^{|\boldsymbol{x}|\sqrt{1-t}}$.

\begin{figure} \caption{Example \ref{sec:parabolic} numerical results in 10 dimensions}\label{fig:5.6}
\minipage{0.45\textwidth}
  \includegraphics[width=\linewidth]{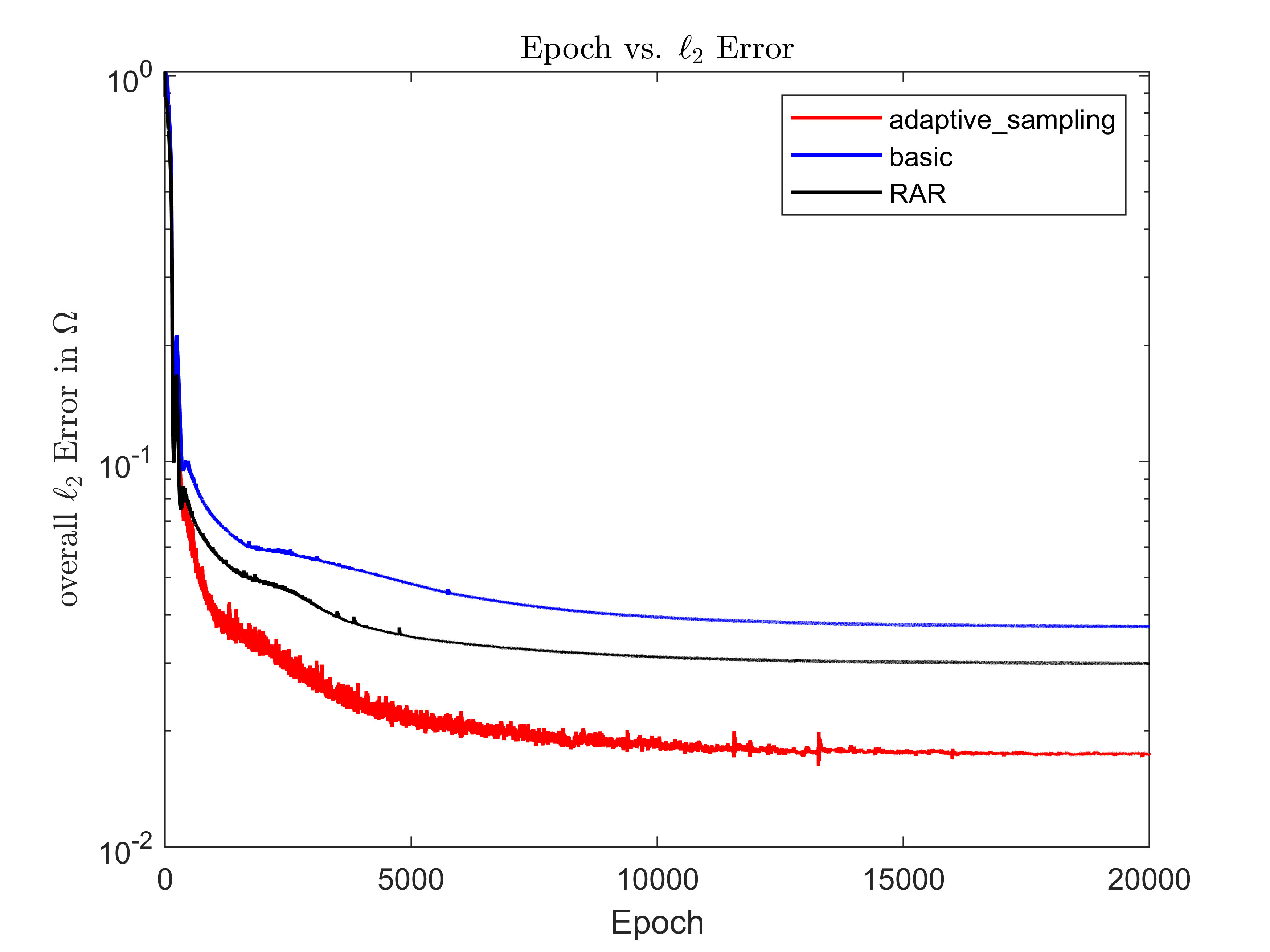}
  \caption*{(a) Epoch vs. Overall Relative $\ell_2$ Error in $\Omega$}
\endminipage\hfil
\minipage{0.45\textwidth}
  \includegraphics[width=\linewidth]{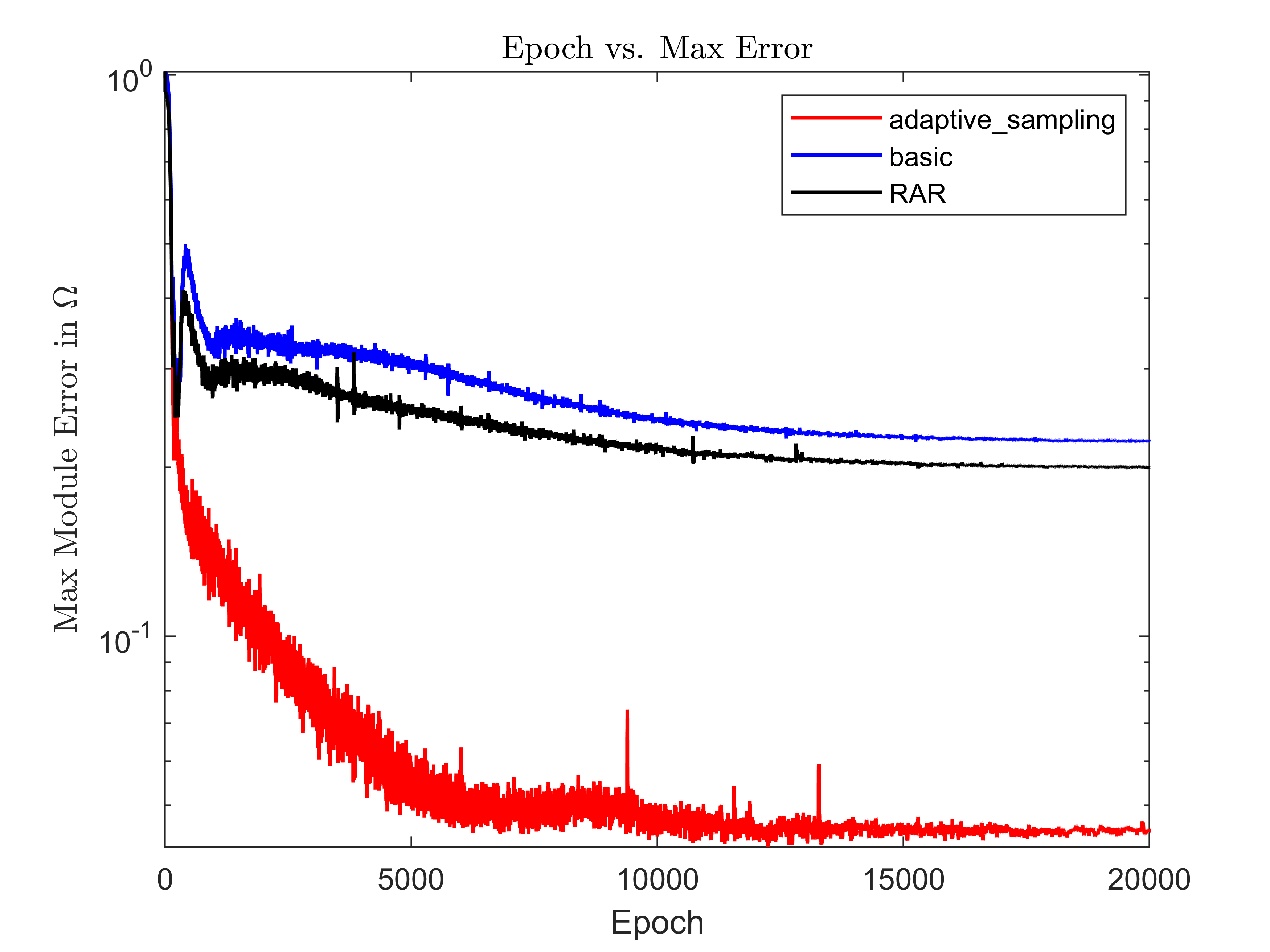}
  \caption*{(b) Epoch vs. Max Relative Modulus Error in $\Omega$}
\endminipage\hfill
\end{figure}

\begin{figure} \caption{Example \ref{sec:parabolic} 10 dimensions $(t, 0, 0, x_3, 0,, ..., 0)$-surface of network solutions and the true solution}\label{fig:5.7}
\minipage{0.33\textwidth}
  \includegraphics[width=\linewidth]{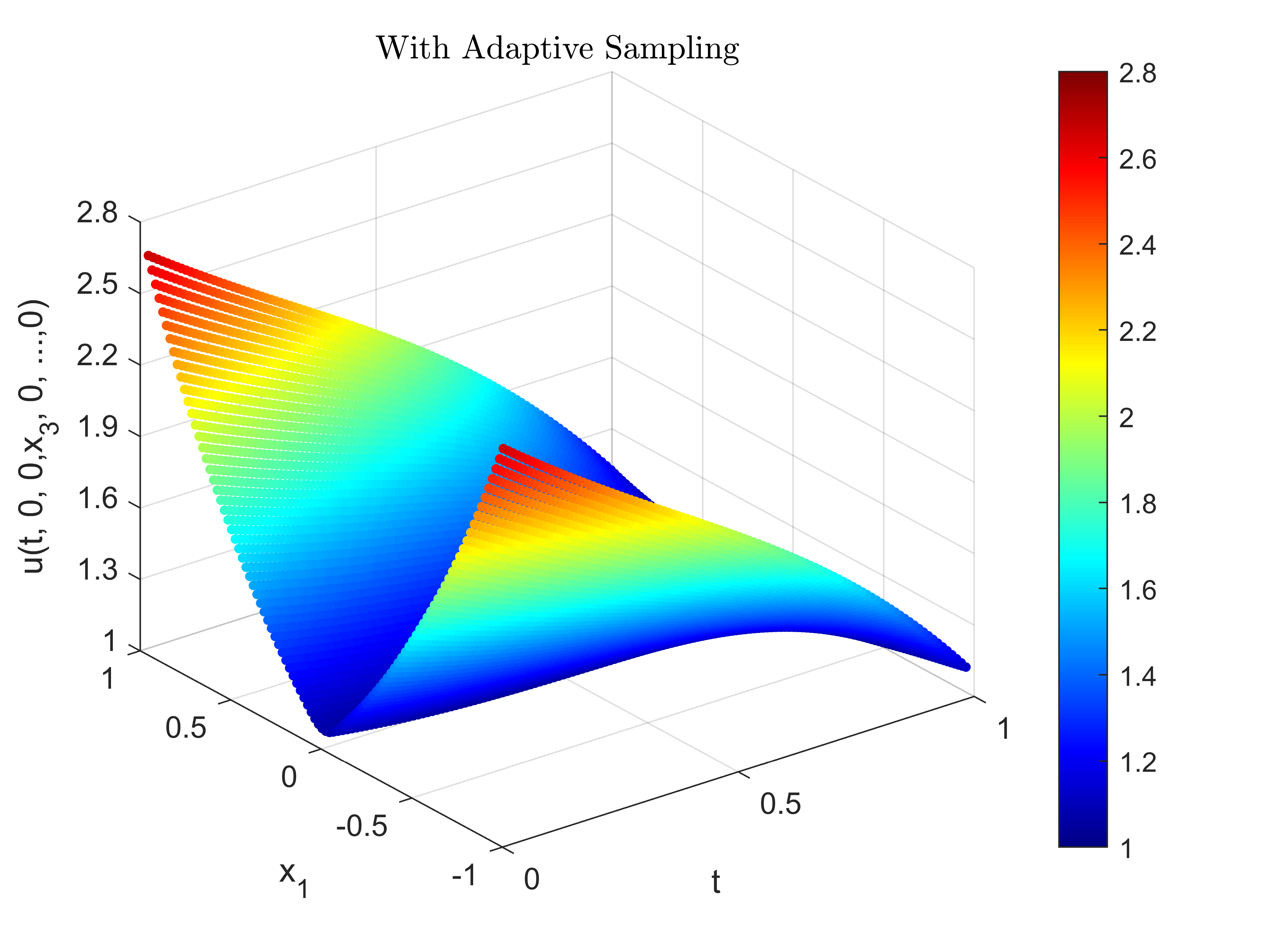}
  \caption*{(a) With Adaptive Sampling}
\endminipage\hfill
\minipage{0.33\textwidth}
  \includegraphics[width=\linewidth]{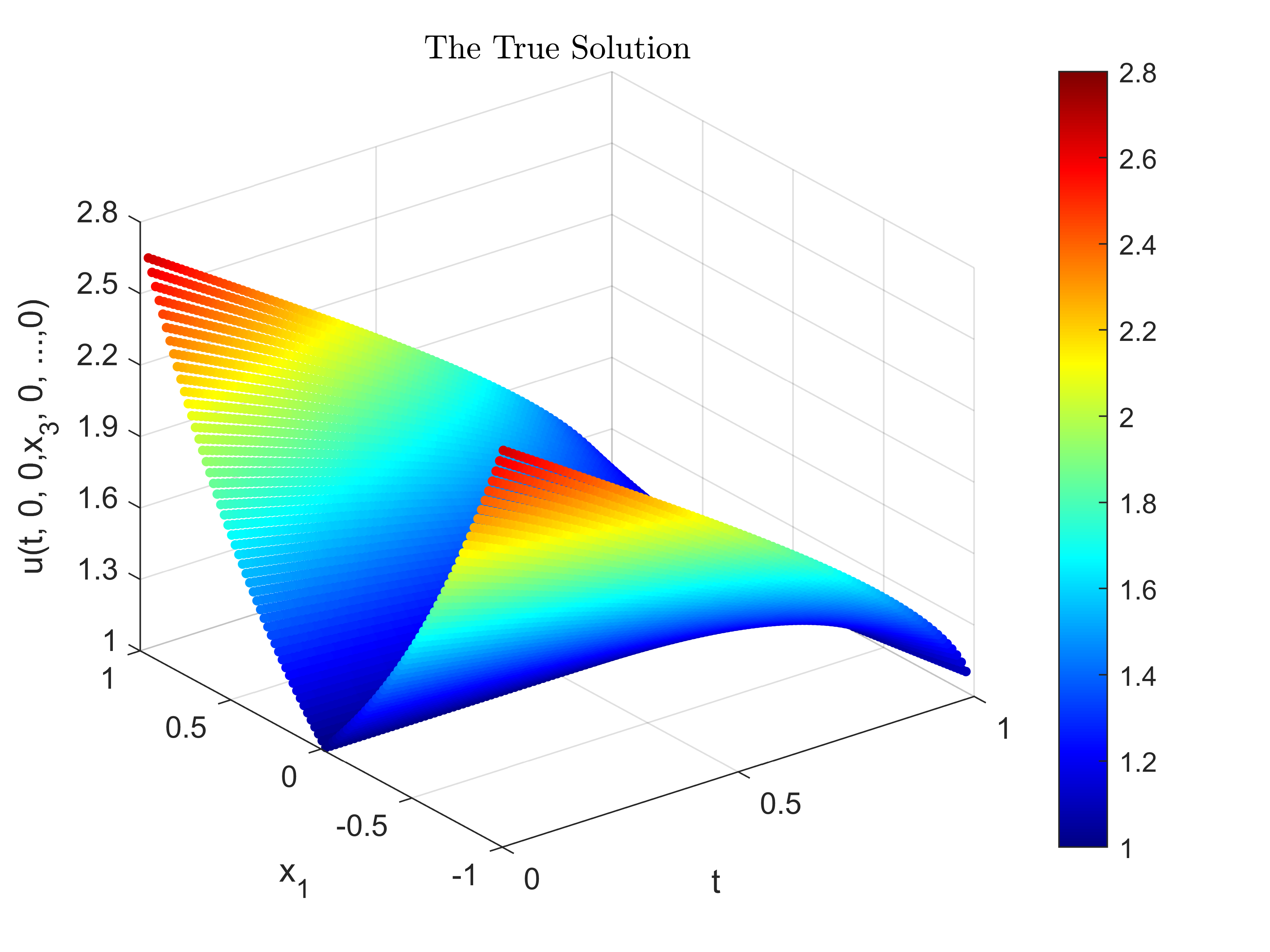}
  \caption*{(b) The True Solution}
\endminipage\hfill
\minipage{0.33\textwidth}
  \includegraphics[width=\linewidth]{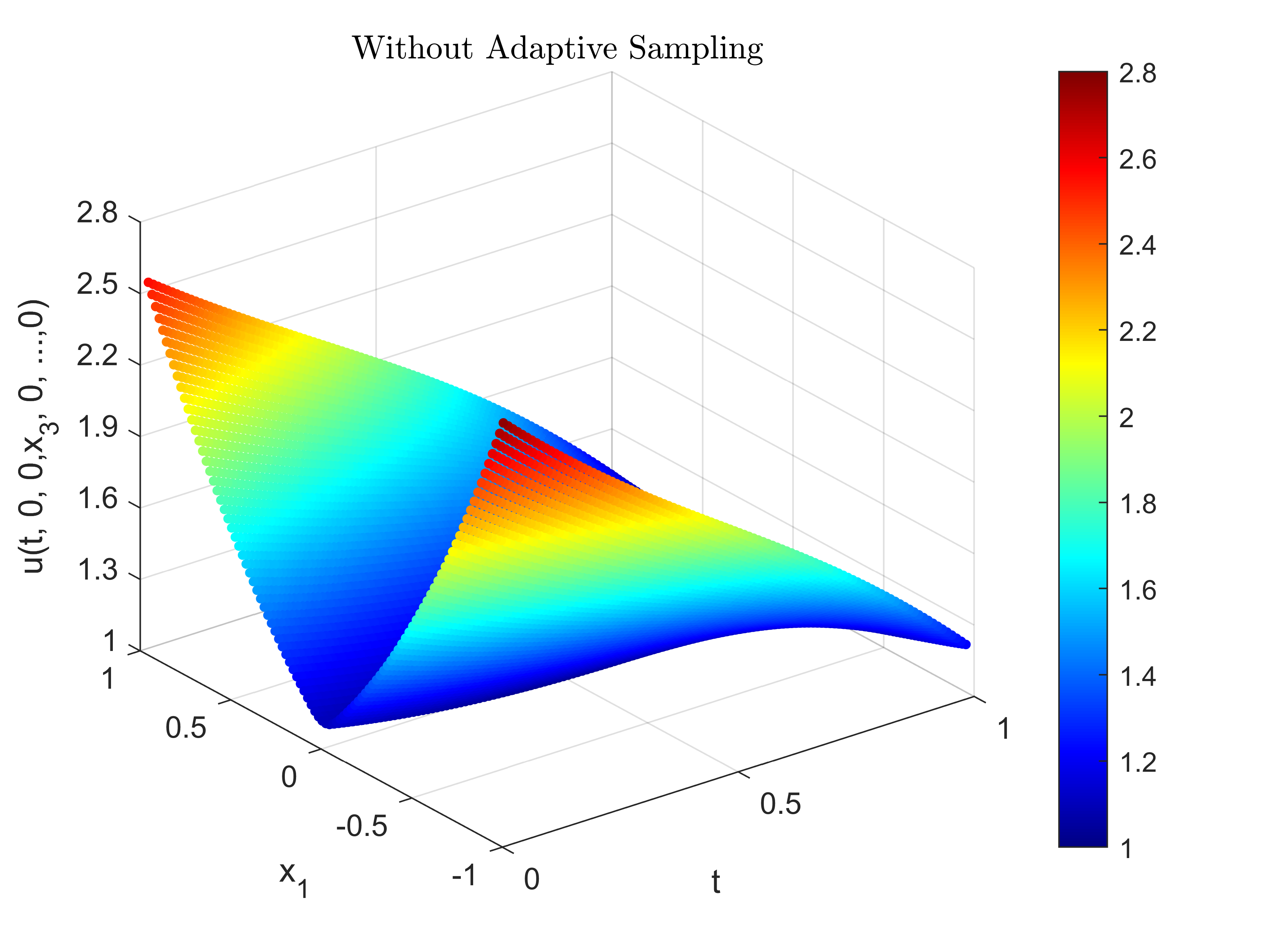}
  \caption*{(c) Without Adaptive Sampling}
\endminipage\hfill
\end{figure}

\begin{figure} \caption{Example \ref{sec:parabolic} 10 dimensions $(t, 0, 0, x_3, 0,, ..., 0)$-surface absolute difference $|u - \phi|$}\label{fig:5.8}
\minipage{0.49\textwidth}
  \includegraphics[width=\linewidth]{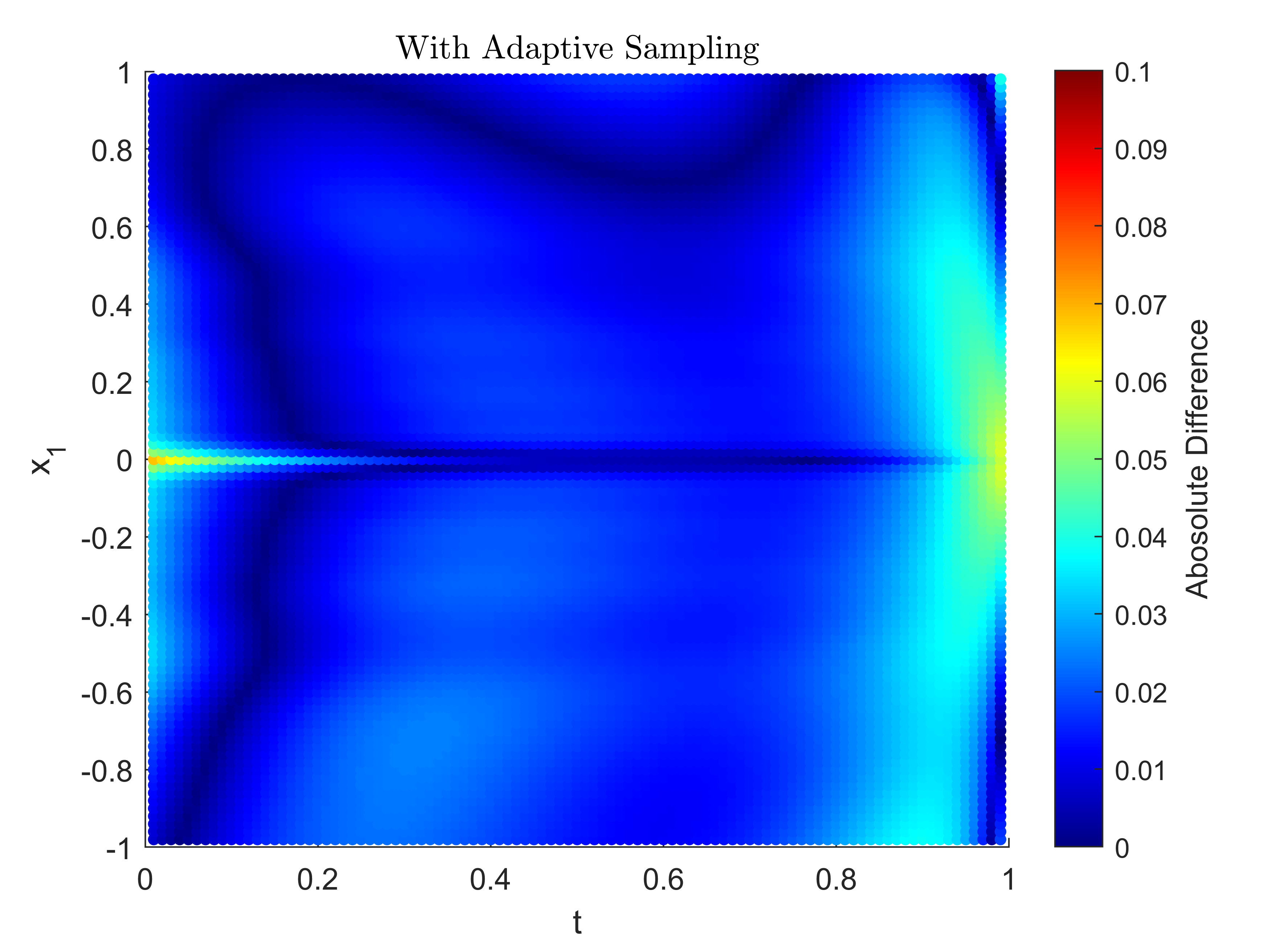}
  \caption*{(a) With Adaptive Sampling}
\endminipage\hfill
\minipage{0.49\textwidth}
  \includegraphics[width=\linewidth]{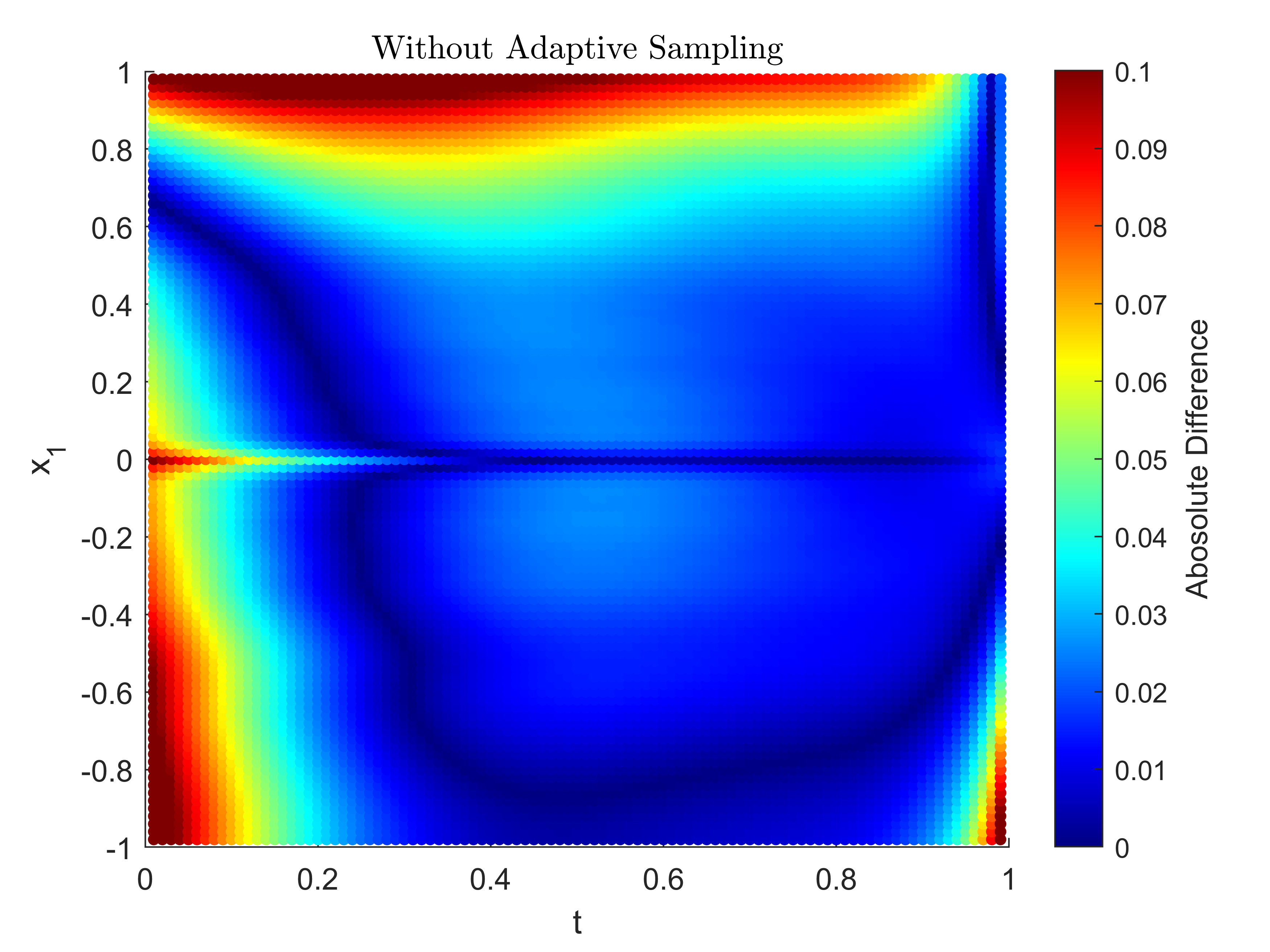}
  \caption*{(b) Without Adaptive Sampling}
\endminipage\hfill
\end{figure}

\begin{figure} \caption{{Example \ref{sec:parabolic} numerical results in 20 dimensions}}\label{fig:5.9}
\minipage{0.45\textwidth}
  \includegraphics[width=\linewidth]{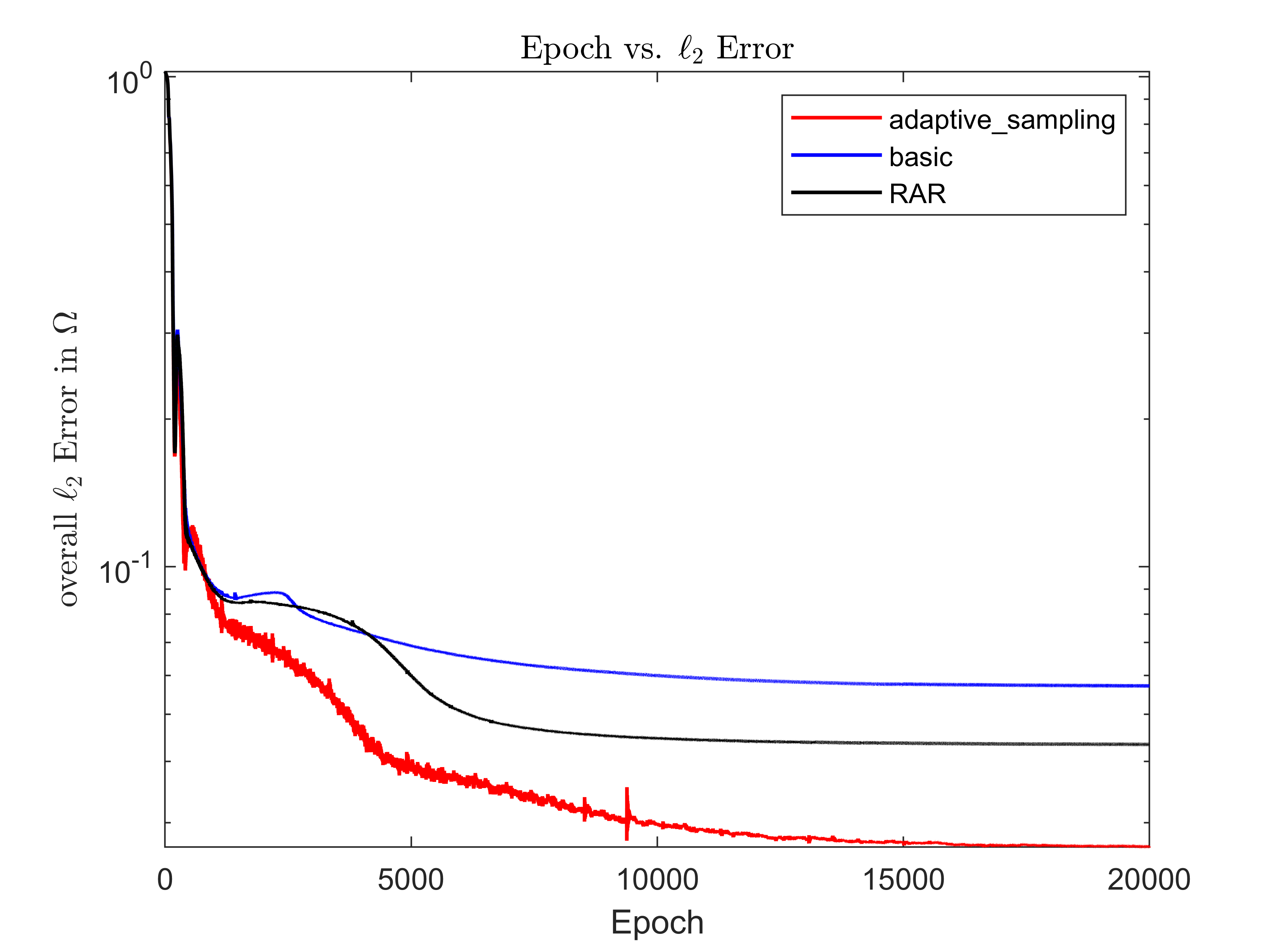}
  \caption*{(a) Epoch vs. Overall Relative $\ell_2$ Error in $\Omega$}
\endminipage\hfil
\minipage{0.45\textwidth}
  \includegraphics[width=\linewidth]{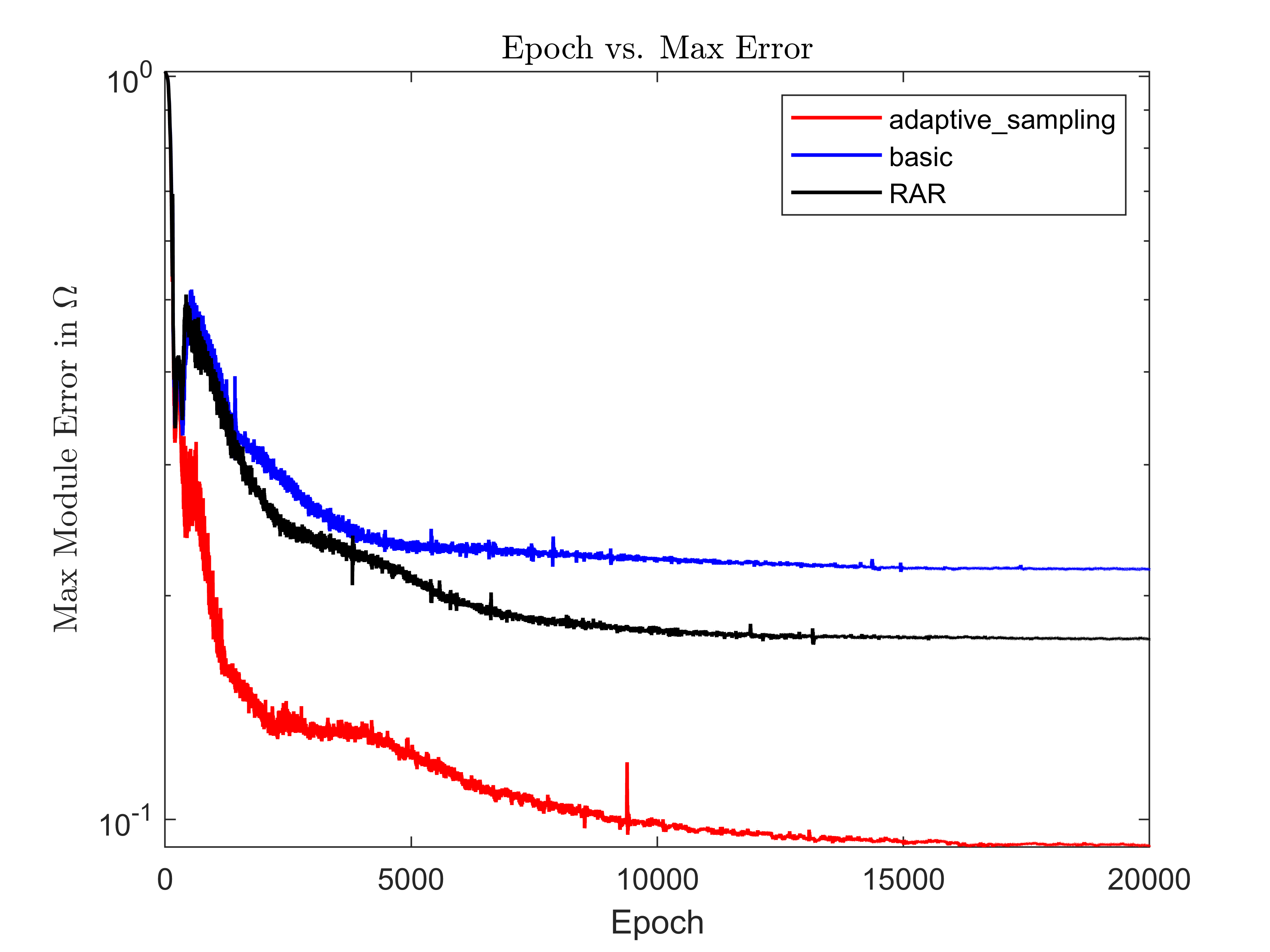}
  \caption*{(b) Epoch vs. Max Relative Modulus Error in $\Omega$}
\endminipage\hfill
\end{figure}

Table \ref{table:5} shows the overall relative $l_2$ error and maximum relative $l_1$ modulus error by adaptive sampling and the  basic residual model  at the end of 20000 epochs. The number of epochs versus $\ell_2$ error decay in 10 dimensions and 20 dimensions  are illustrated in Figures \ref{fig:5.6} and \ref{fig:5.9} respectively. 
Figure \ref{fig:5.7} shows the $(t, x_3)$-surface of the ground-truth solution and network solutions with and without adaptive sampling, where axes are $t$, $x_3$, and $u(t, 0, 0, x_3, 0, ..., 0)$. Figure \ref{fig:5.8} shows the heat-map of the absolute difference in $(t, x_3)$-surface, where the color bar represents the absolute difference between the ground-truth solution and network solution. This example evidently demonstrates the advantage of adaptive sampling in reducing the variance in the distribution of error over the domain. These results clearly show that adaptive sampling helps to reduce the error, especially relative $l_1$ maximum modulus error.

\begin{table}
	\caption{Example \ref{sec:parabolic} Result}\label{table:5}
	\centering
	\begin{tabular}{llllll}
		\toprule
		Dimension& & ~~ AS & ~~Basic & ~~RAR & ~~Error Reduced by AS\\
		\midrule
		 10 & \begin{tabular}{@{}l@{}}
                   $\ell_2$ error\\
                   Max Modulus Error\\
                   \\
                 \end{tabular}& \begin{tabular}{l@{}@{}}
                   1.731916e-02\\
                   4.530897e-02\\
                   \\
                 \end{tabular}  &\begin{tabular}{l@{}@{}}
                   3.735915e-02\\
                   2.229580e-01\\
                   \\
                 \end{tabular}   &\begin{tabular}{l@{}@{}}
                   2.997037e-02\\
                   2.000248e-01\\
                   \\
                 \end{tabular}     
                 &\begin{tabular}{l@{}@{}}
                   53.54$\%$\\
                   79.68$\%$\\
                   \\
                 \end{tabular}     \\
		
		 20 & \begin{tabular}{@{}l@{}}
                   $\ell_2$ error\\
                   Max Modulus Error\\
                   \\
                 \end{tabular}& \begin{tabular}{l@{}@{}}
                   2.688834e-02\\
                   9.229024e-02\\
                   \\
                 \end{tabular}  &\begin{tabular}{l@{}@{}}
                   5.709440e-02\\
                   2.167297e-01\\
                   \\
                 \end{tabular}  &\begin{tabular}{l@{}@{}}
                   4.336126e-02\\
                   1.748666e-01\\
                   \\
                 \end{tabular}  &\begin{tabular}{l@{}@{}}
                   52.91$\%$\\
                   57.42$\%$\\
                   \\
                 \end{tabular}     \\
		\bottomrule
	\end{tabular}
\end{table}

\subsection{Hyperbolic Equation}\label{sec: hyperb}
\noindent Lastly, we consider the following hyperbolic equation:
\begin{equation} \label{eq:hyperbo}
\begin{aligned}
\frac{\partial^{2} u(x, t)}{\partial t^{2}}-\Delta_{x} u(x, t)=f(x, t)&,\quad \text { in } \Omega:=\omega \times \mathbb{T}, \\
u(x, t)=g_{0}(x, t)&,\quad \text { on } \partial \Omega = \partial \omega \times \mathbb{T}, \\
u(x, 0)=h_{0}(x) \quad \frac{\partial u(x, 0)}{\partial t}=h_{1}(x)&,\quad \text { in } \omega,
\end{aligned}
\end{equation}
 where $\omega:=\{x:|x|<1\}$, $\mathbb{T} = (0, 1)$. 
 $g_{0}(x, t) = 0$, $h_{0}(x) = 0$, $h_{1}(x) = 0$, 
and $f(x,t)$ is given appropriately so that 
the exact solution is
 $u(x, t)=\left(\exp \left(t^{2}\right)-1\right) \sin \left(\frac{\pi}{2}(1-|x|)^{2.5}\right)
$.  $\Delta_{x}$ denotes the Laplace operator taken in the spatial variable $x$ only.

\begin{figure} \caption{Example \ref{sec: hyperb} numerical results in 10 dimensions}\label{fig:5.10}
\minipage{0.45\textwidth}
  \includegraphics[width=\linewidth]{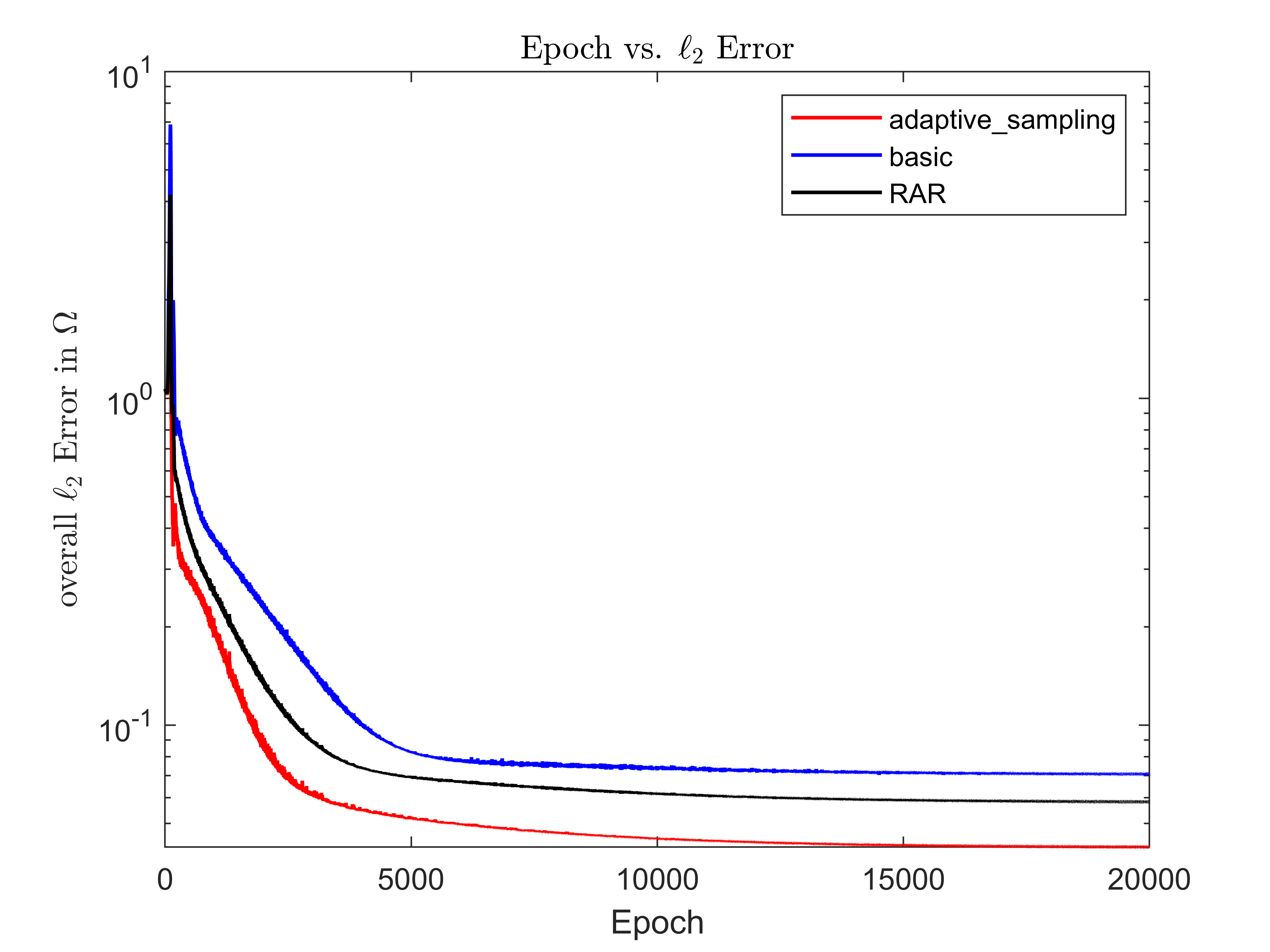}
  \caption*{(a) Epoch vs. Overall Relative $\ell_2$ Error in $\Omega$}
\endminipage\hfil
\minipage{0.45\textwidth}
  \includegraphics[width=\linewidth]{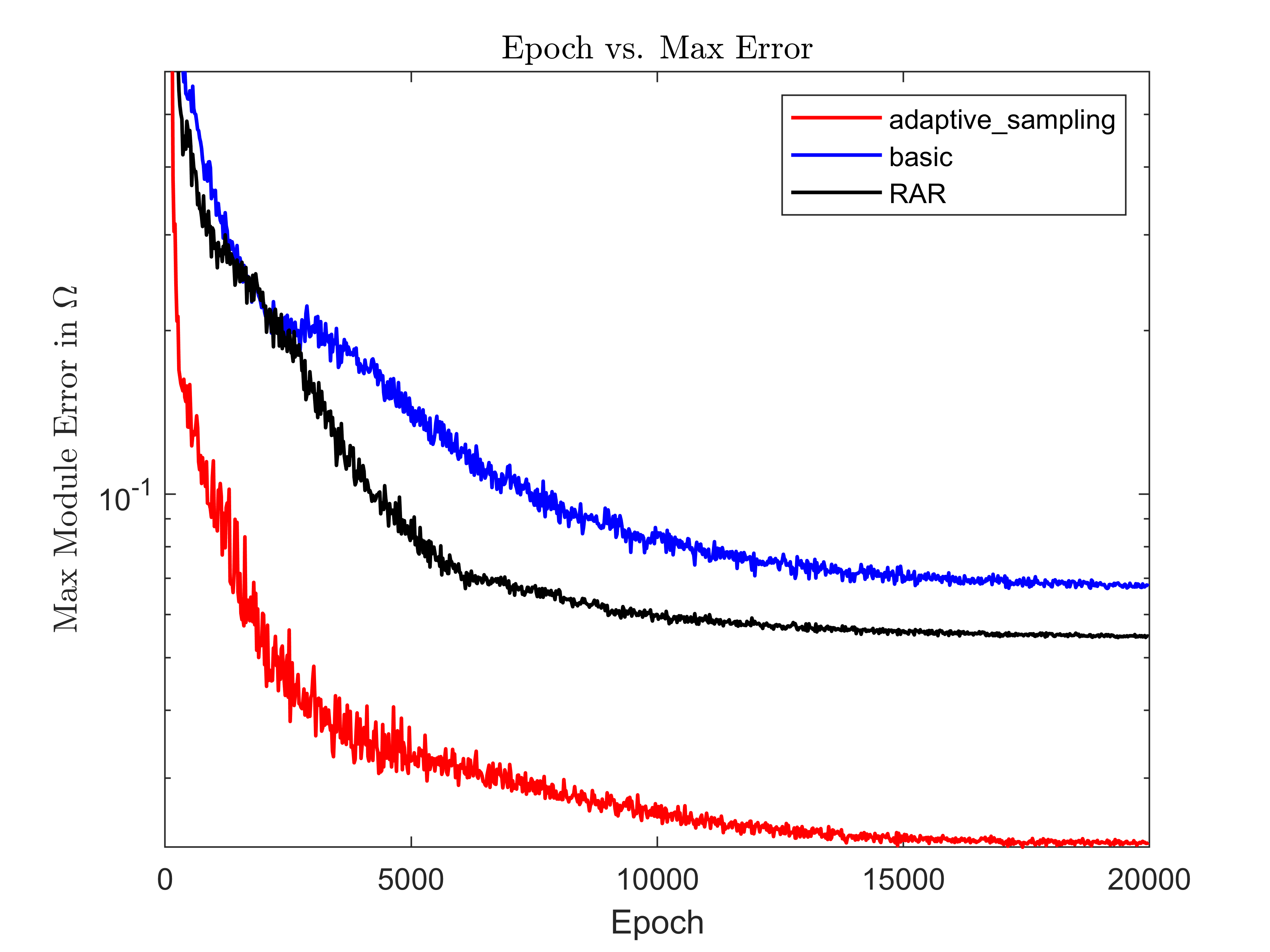}
  \caption*{(b) Epoch vs. Max Relative Modulus Error in $\Omega$}
\endminipage\hfill
\end{figure}

\begin{figure} \caption{Example \ref{sec: hyperb} 10 dimensions $(t, 0,...,0, x_6,0,...,0)$-surface of network solutions and the true solution}\label{fig:5.11}
\minipage{0.33\textwidth}
  \includegraphics[width=\linewidth]{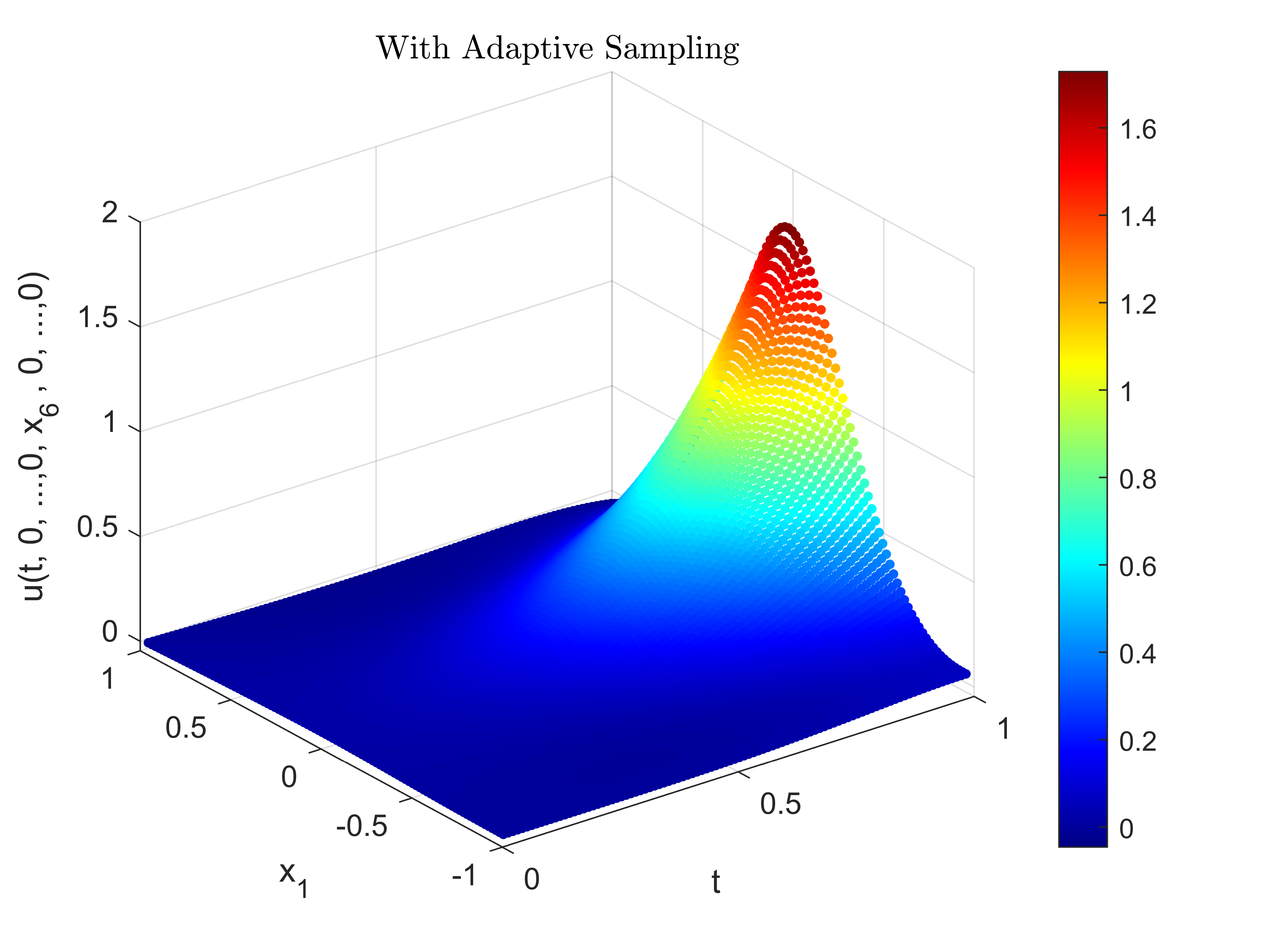}
  \caption*{(a) With Adaptive Sampling}
\endminipage\hfill
\minipage{0.33\textwidth}
  \includegraphics[width=\linewidth]{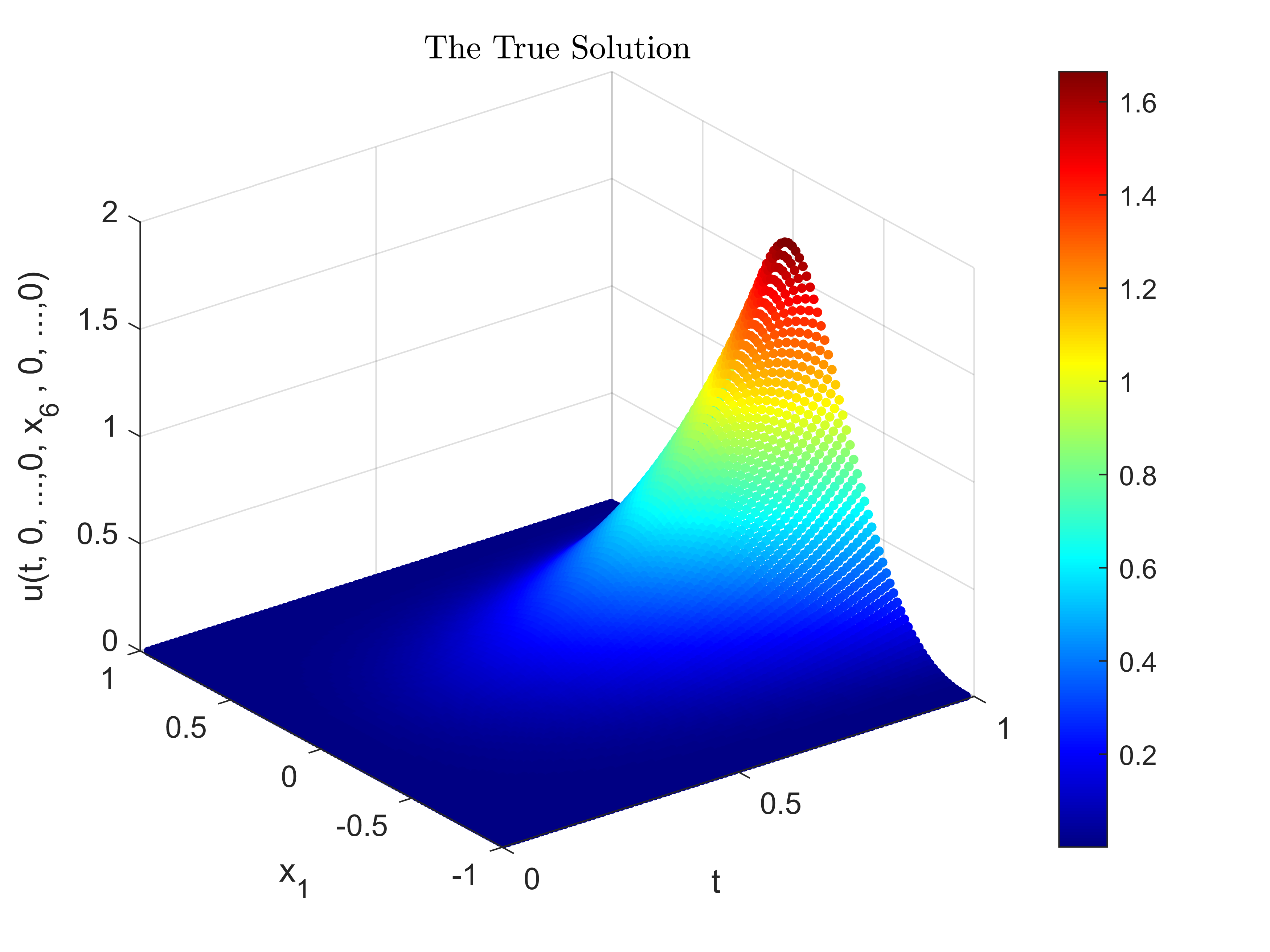}
  \caption*{(b) The True Solution}
\endminipage\hfill
\minipage{0.33\textwidth}
  \includegraphics[width=\linewidth]{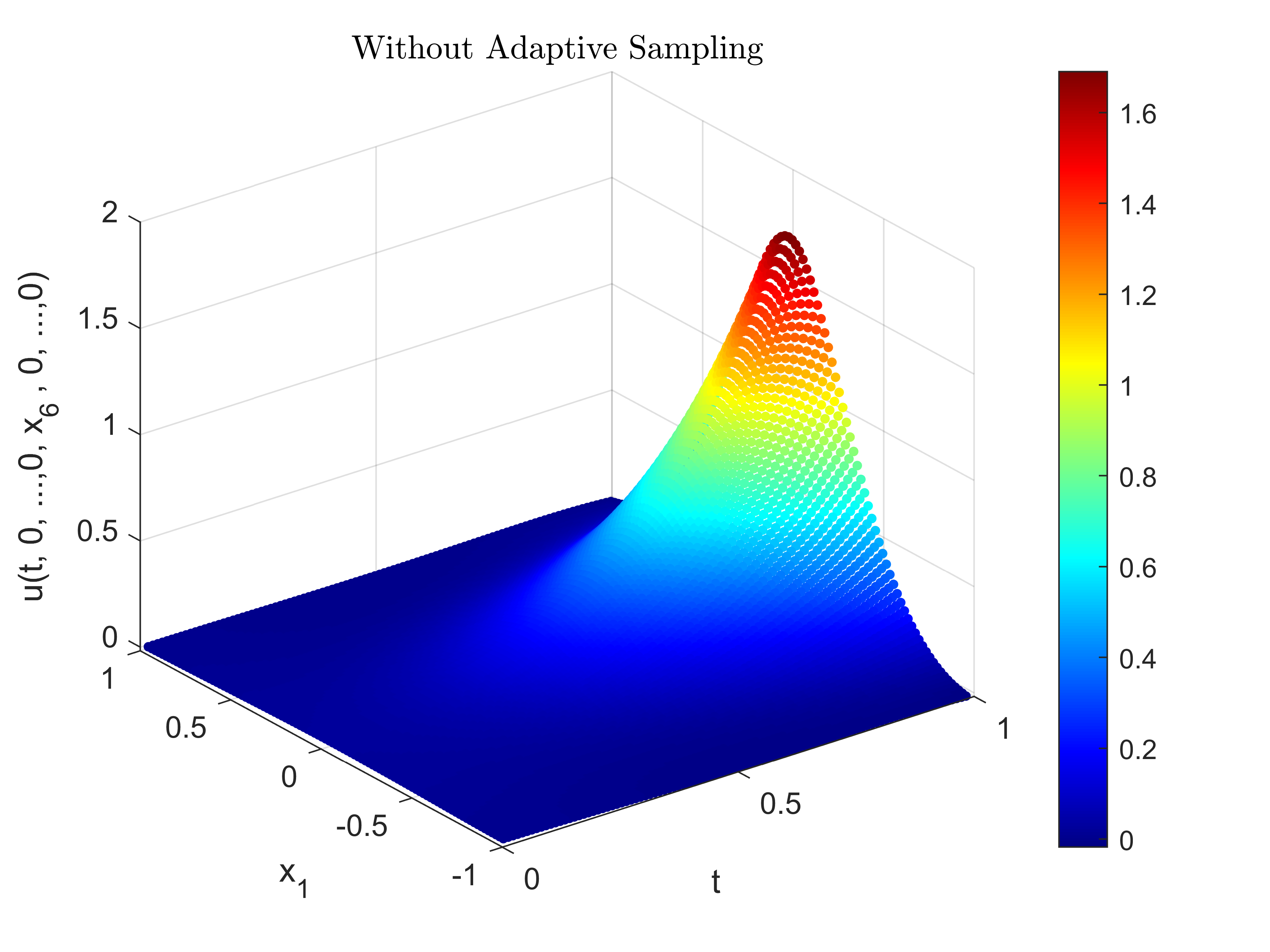}
  \caption*{(c) Without Adaptive Sampling}
\endminipage\hfill
\end{figure}

\begin{figure} \caption{Example \ref{sec: hyperb} 10 dimensions $(t, 0,..., 0, x_6, 0,, ..., 0)$-surface absolute difference $|u - \phi|$}\label{fig:5.12}
\minipage{0.49\textwidth}
  \includegraphics[width=\linewidth]{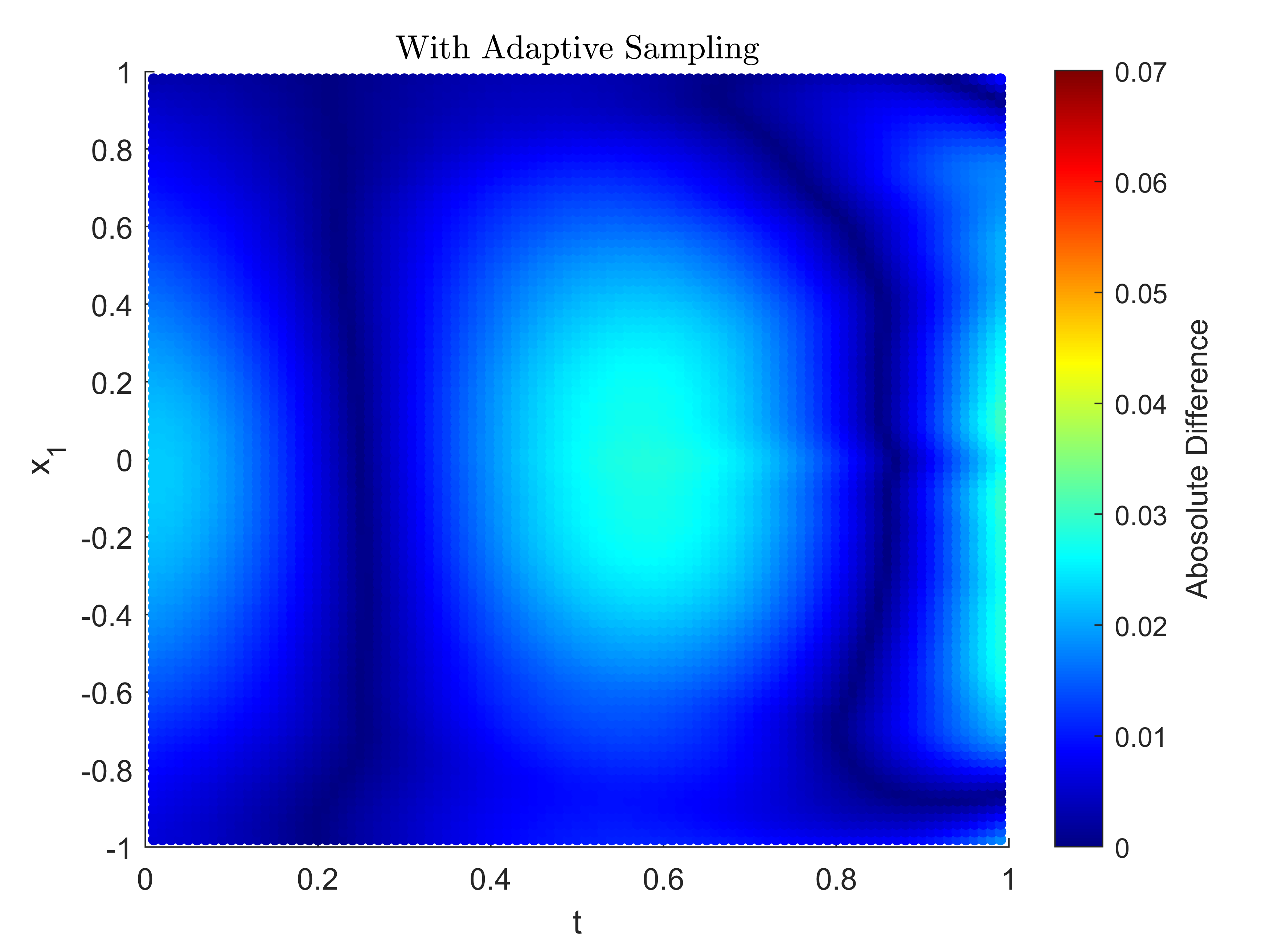}
  \caption*{(a) With Adaptive Sampling}
\endminipage\hfill
\minipage{0.49\textwidth}
  \includegraphics[width=\linewidth]{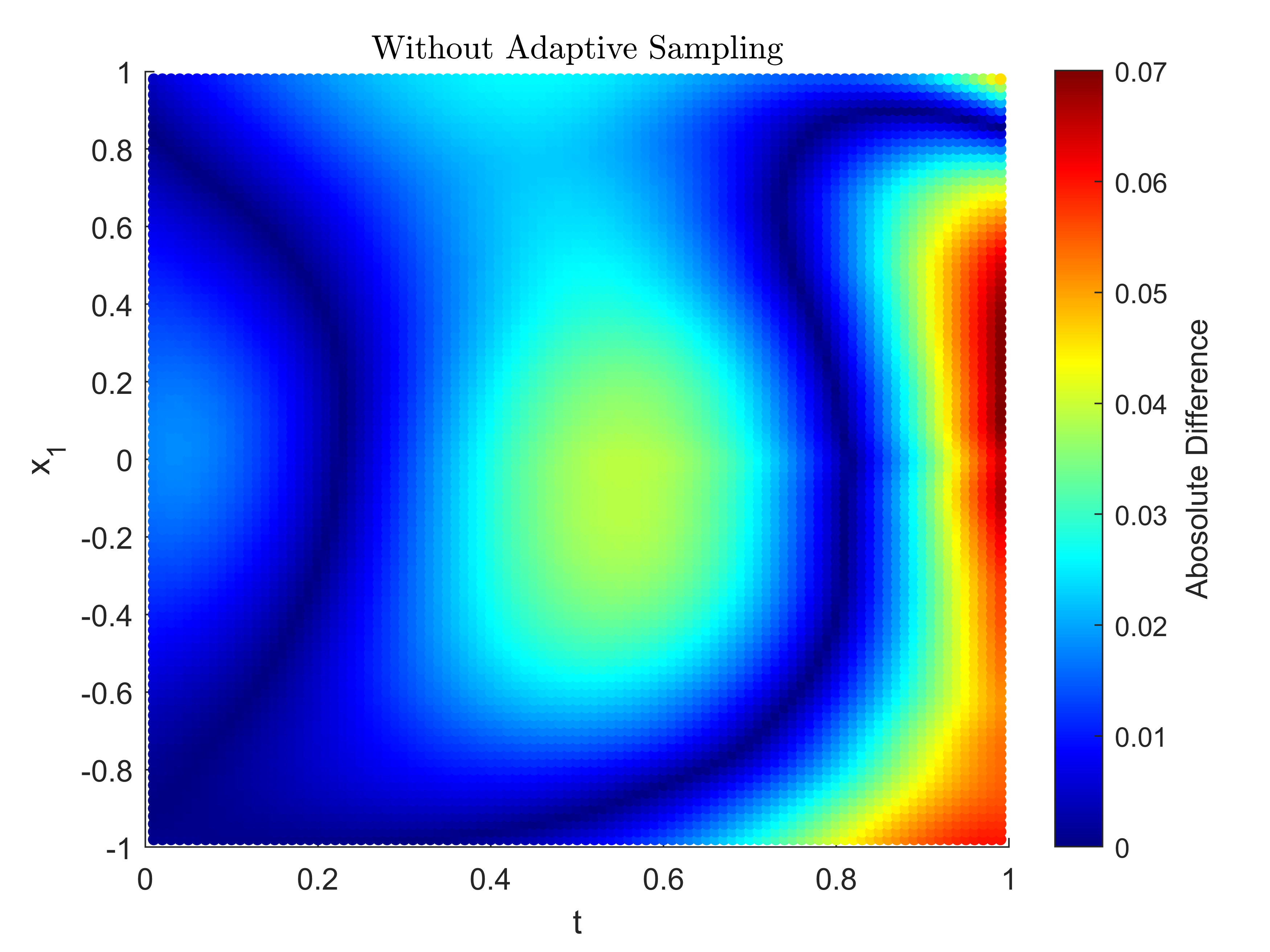}
  \caption*{(b) Without Adaptive Sampling}
\endminipage\hfill
\end{figure}

\begin{figure} \caption{Example \ref{sec: hyperb} numerical results in 20 dimensions}\label{fig:5.13}
\minipage{0.45\textwidth}
  \includegraphics[width=\linewidth]{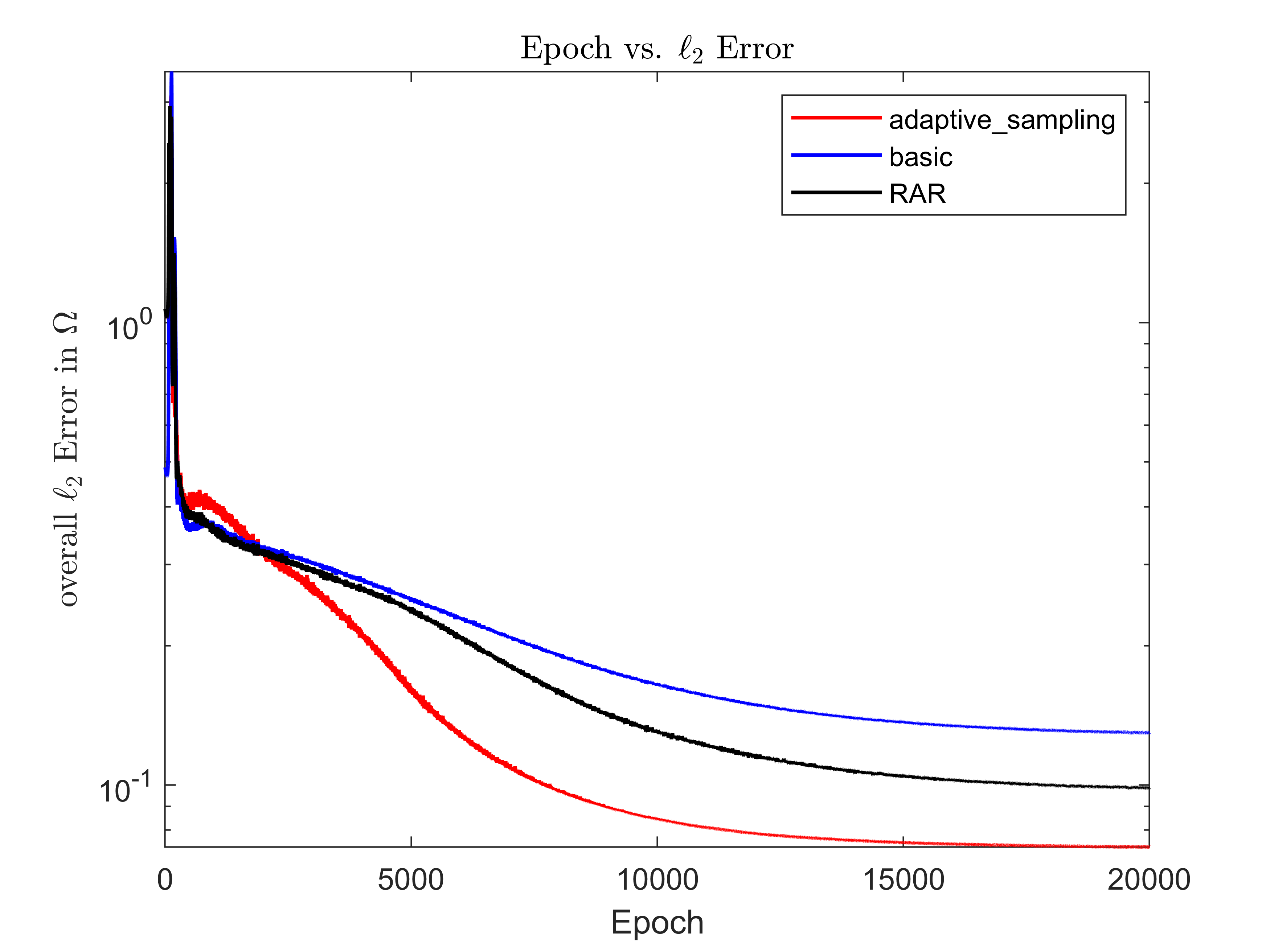}
  \caption*{(a) Epoch vs. Overall Relative $\ell_2$ Error in $\Omega$}
\endminipage\hfil
\minipage{0.45\textwidth}
  \includegraphics[width=\linewidth]{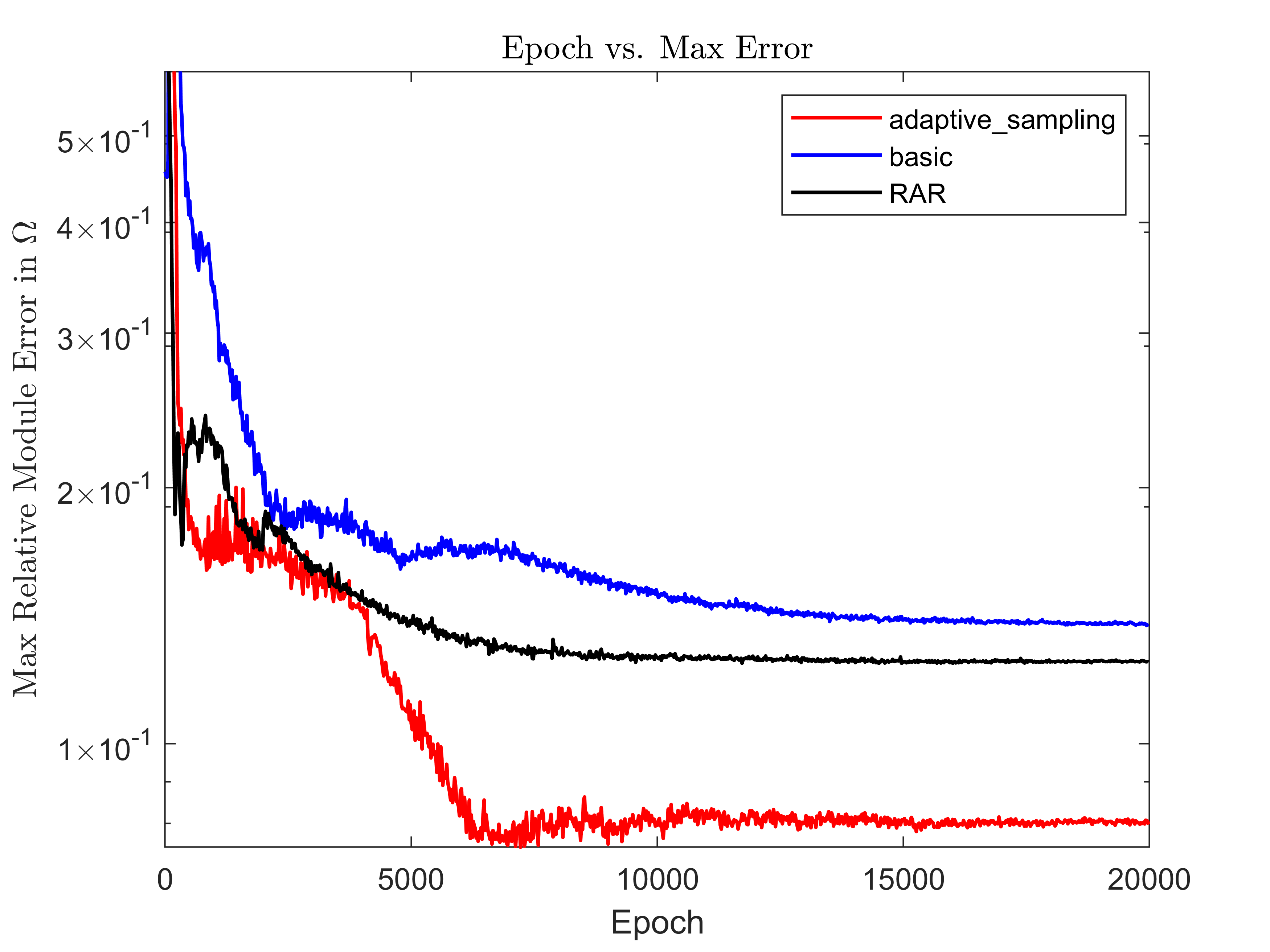}
  \caption*{(b) Epoch vs. Max Relative Modulus Error in $\Omega$}
\endminipage\hfill
\end{figure}

Table \ref{table:6} shows the overall relative $l_2$ error and maximum relative $l_1$ modulus error by adaptive sampling and the basic model at the end of 20000 epochs.  The number of epochs versus $\ell_2$ error decay in 10 dimensions and 20 dimensions  are presented respectively in  Figure \ref{fig:5.10} and \ref{fig:5.13}. 
Figure \ref{fig:5.11} shows the $(t, x_6)$-surface of the ground-truth solution and network solutions with and without adaptive sampling, where axes are $t$, $x_6$, and $u(t, 0,..., 0, x_6, 0,, ..., 0)$. Figure \ref{fig:5.12} shows the heat-map of absolute difference in $(t, x_6)$-surface, where the color bar represents the absolute difference between the ground-truth solution and network solution. This example  shows that the adaptive sampling  reduces the variance in the distribution of error over the domain. We can see from these results that adaptive sampling helps to reduce the error, especially relative $l_1$ maximum modulus error.

\begin{table}
	\caption{Example \ref{sec: hyperb} Result}\label{table:6}
	\centering
	\begin{tabular}{llllll}
		\toprule
		Dimension& &~~AS& ~~Basic & ~~RAR & ~~Error Reduced by AS\\
		\midrule
		 10 & \begin{tabular}{@{}l@{}}
                   $\ell_2$ error\\
                   Max Modulus Error\\
                   \\
                 \end{tabular}& \begin{tabular}{l@{}@{}}
                   4.233210e-02\\
                   2.265792e-02\\
                   \\
                 \end{tabular}  &\begin{tabular}{l@{}@{}}
                   7.093248e-02\\
                   6.659816e-02\\
                   \\
                 \end{tabular}   
                 &\begin{tabular}{l@{}@{}}
                   5.827290e-02\\
                   5.484901e-02'\\
                   \\
                 \end{tabular}
                 &\begin{tabular}{l@{}@{}}
                   40.32$\%$\\
                   65.98$\%$\\
                   \\
                 \end{tabular}     \\
		
		 20 & \begin{tabular}{@{}l@{}}
                   $\ell_2$ error\\
                   Max Modulus Error\\
                   \\
                 \end{tabular}& \begin{tabular}{l@{}@{}}
                   7.342700e-02\\
                   8.995471e-02\\
                   \\
                 \end{tabular}  &\begin{tabular}{l@{}@{}}
                   1.299070e-01\\
                   1.484266e-01\\
                   \\
                 \end{tabular} 
                 &\begin{tabular}{l@{}@{}}
                   9.876399e-02\\
                  1.350940e-01\\
                   \\
                 \end{tabular}    
                 &\begin{tabular}{l@{}@{}}
                   43.48$\%$\\
                   39.39$\%$\\
                   \\
                 \end{tabular}     \\
                 
		\bottomrule
	\end{tabular}
\end{table}

\subsection{Discussion on Comparison with RAR}
In all three PDE examples above, AS and RAR both can speed up the convergence and lower the error. In the above results, AS performs better than RAR in high-dimensional settings. Here we provide some intuition behind this out-performance. As we mentioned in Section \ref{sec:relevant lit.}, AS differentiates all training points and is more prone to select points with higher residual errors. On the other hand, RAR does not differentiate between training points; RAR replicates $k$ points with the largest residuals. In AS, by ranking all training points, the neural network model focuses more on points that the model is more uncertain about. As we can see from the above results, even though AS does a good job in reducing the overall relative $\ell_2$ error, AS is even more effective in reducing the max relative modulus error. This is because when the model is producing high residual errors on some volumes, AS will make the model highly focused on these volumes. In comparison, although RAR also focuses on points with high residuals, the effect might not be enough, especially in high-dimensional cases. As mentioned in Section \ref{sec:overview_of_AL}, there are two main approaches to active learning: uncertainty sampling and diversity sampling. AS is an uncertainty sampling approach. Intuitively, RAR is more diversity-preserving than AS because it still keeps points that are randomly generated. However, AS also preserves diversity throughout the training process. 
AS will tackle points with higher residuals first, but after some iterations, the residual error of these points will be reduced, and other points that were previously not considered ``high-residual" points now will be considered ``high-residual" points by AS. Therefore, over the course of the whole training process, the AS algorithm will balance itself in selecting points that not only focuses on the uncertainty of the neural network model but also preserves diversity in some sense. In conclusion, AS outperforms RAR because AS focuses more on ``high-residual error" areas/volumes as demonstrated in the massive reduction of max modulus error. AS also balances itself in selecting training points over the course of the whole training process so that AS will still preserve diversity to some extent.

\subsection{Poisson's Equation} \label{sec:5.4}
In this example, we will 
show that adaptive sampling is compatible with some recent frameworks: the Deep Ritz Method (DRM)\cite{Deep_Ritz}, the Deep Galerkin Method (DGM)\cite{other3}, and the Weak Adversarial Networks (WAN)\cite{WAN}. We still use residual errors to define the sampling distribution and use Algorithm \ref{alg:Self-normalized Sampling} to sample training points for all frameworks. To this end, we consider a 2D Poisson's equation: find $u(x, y)$ such that
\begin{eqnarray*}
    -\Delta u(x, y) = f(x, y),& \quad& \text { in } \Omega:=(0, 1) \times (0, 1), \\
    u(x, y) = g(x, y),& \quad&  \text { on } \partial \Omega. 
\end{eqnarray*}
 The exact solution is given by $u(x, y) = \min\{x^2, (1-x)^2\}$. All models are trained for 300 seconds in NVIDIA Quadro P1000 Graphics Card. Each framework  without and with adaptive sampling both have 30 trials to approximate the PDE.
The error obtained with adaptive sampling is compared with the error obtained without adaptive sampling to test if adaptive sampling is compatible with various frameworks to lower the error. Note that, in different trials, different random seeds are used. Moreover,  in the $k$-th trial, the same torch random seed is used  with and without adaptive sampling, which implies that the networks are initialized with the same parameters.
In this test, ResNet  \cite{he2015deep} with 4 residual blocks is used with Tanh as activation functions. Each epoch has 1024 uniform points on the boundary and in the domain respectively for training. Moreover, instead of numerical differentiation, Pytorch Autograd is utilized to evaluate derivatives.

\begin{figure}  \caption{Compatibility Test: Five Lowest $\ell_2$ Error Trials in Each Framework}\label{fig:5.111}
\minipage{0.49\textwidth}
  \includegraphics[width=\linewidth]{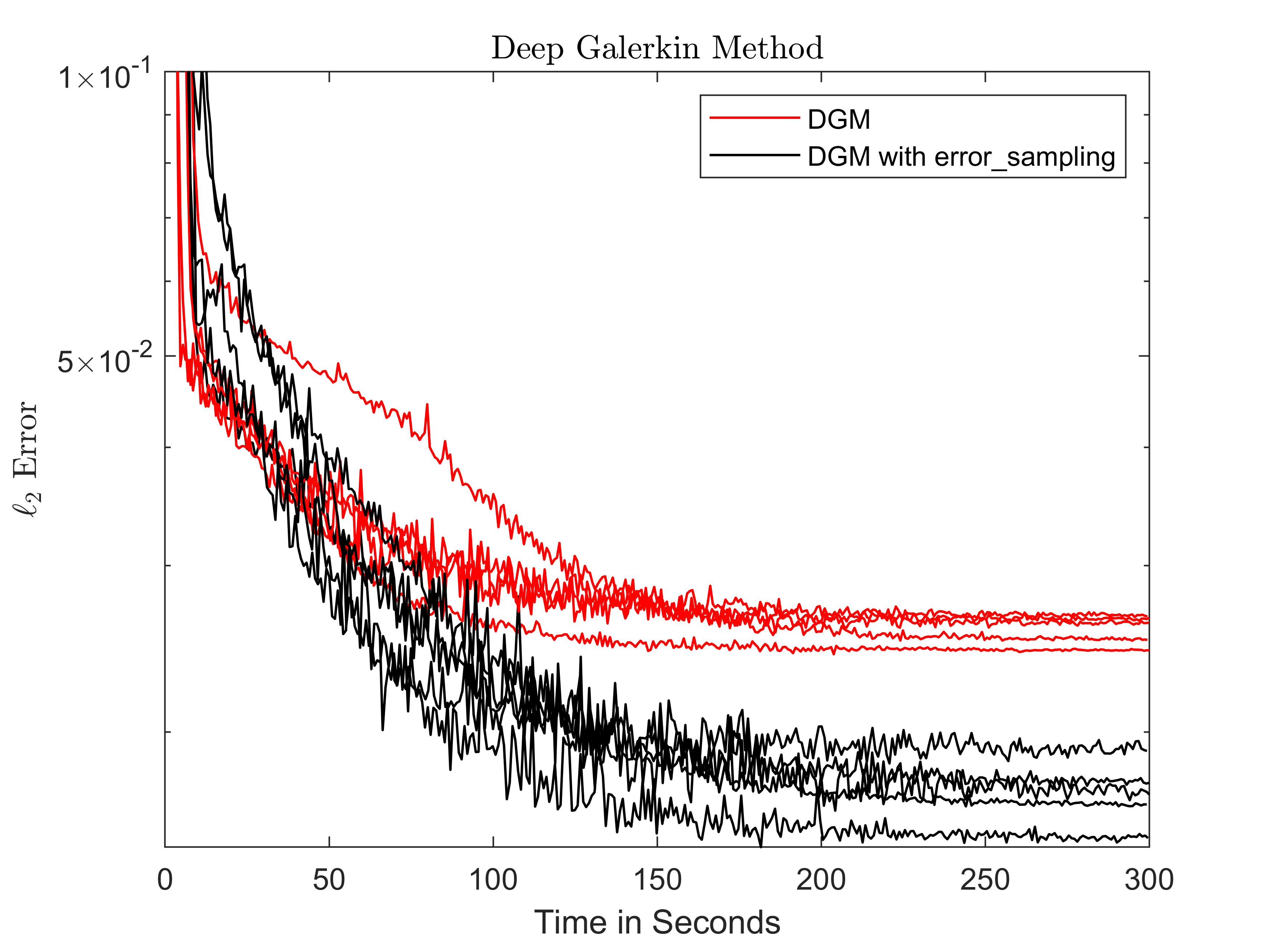}
    \caption*{(a) {DGM: $\ell_2$ Error vs Time}}
\endminipage\hfill
\minipage{0.49\textwidth}
  \includegraphics[width=\linewidth]{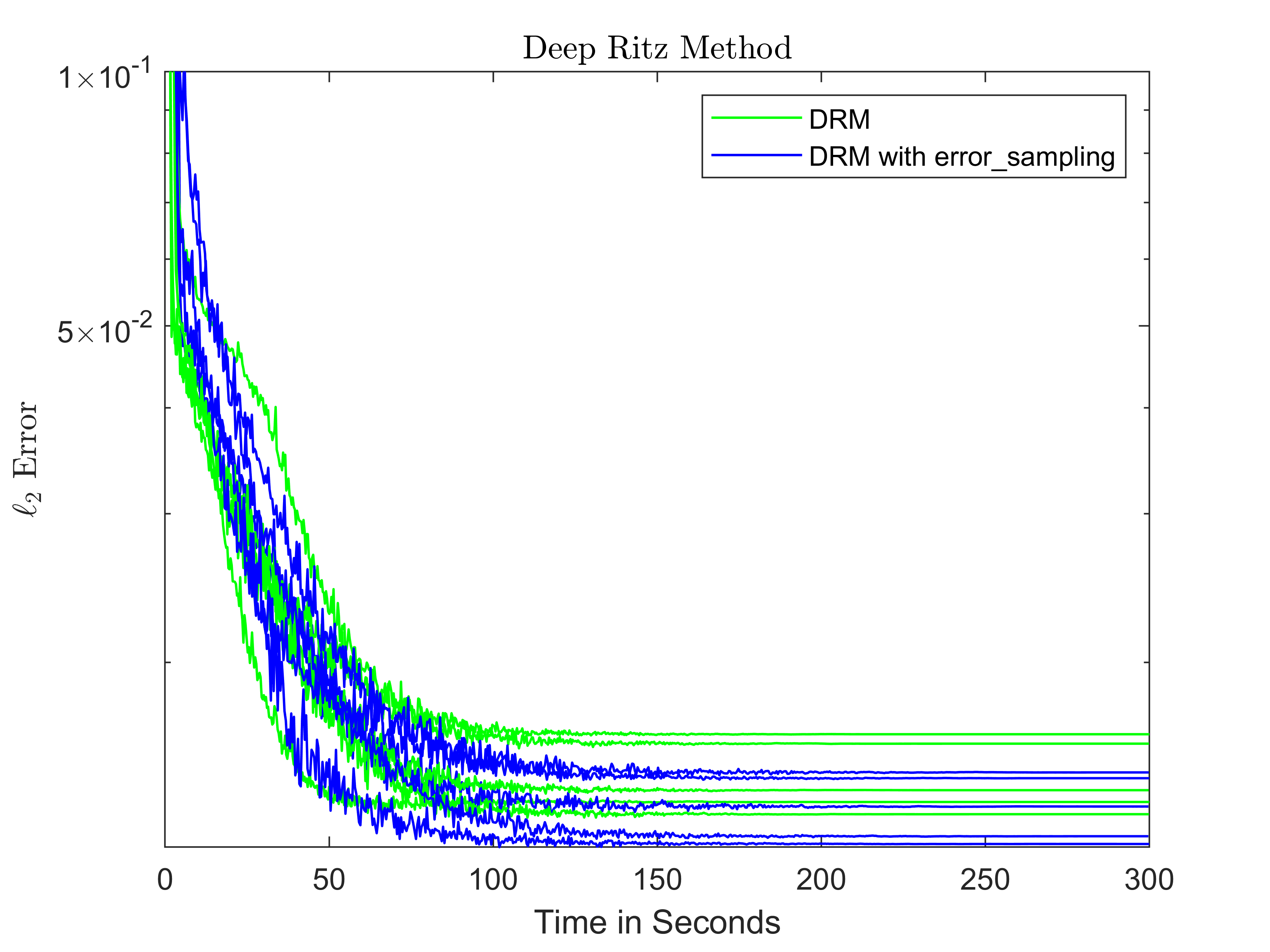}
  \caption*{(b){DRM: $\ell_2$ Error vs Time}}
\endminipage\hfill
\minipage{\textwidth}%
  \centering
  \includegraphics[width=0.49\linewidth]{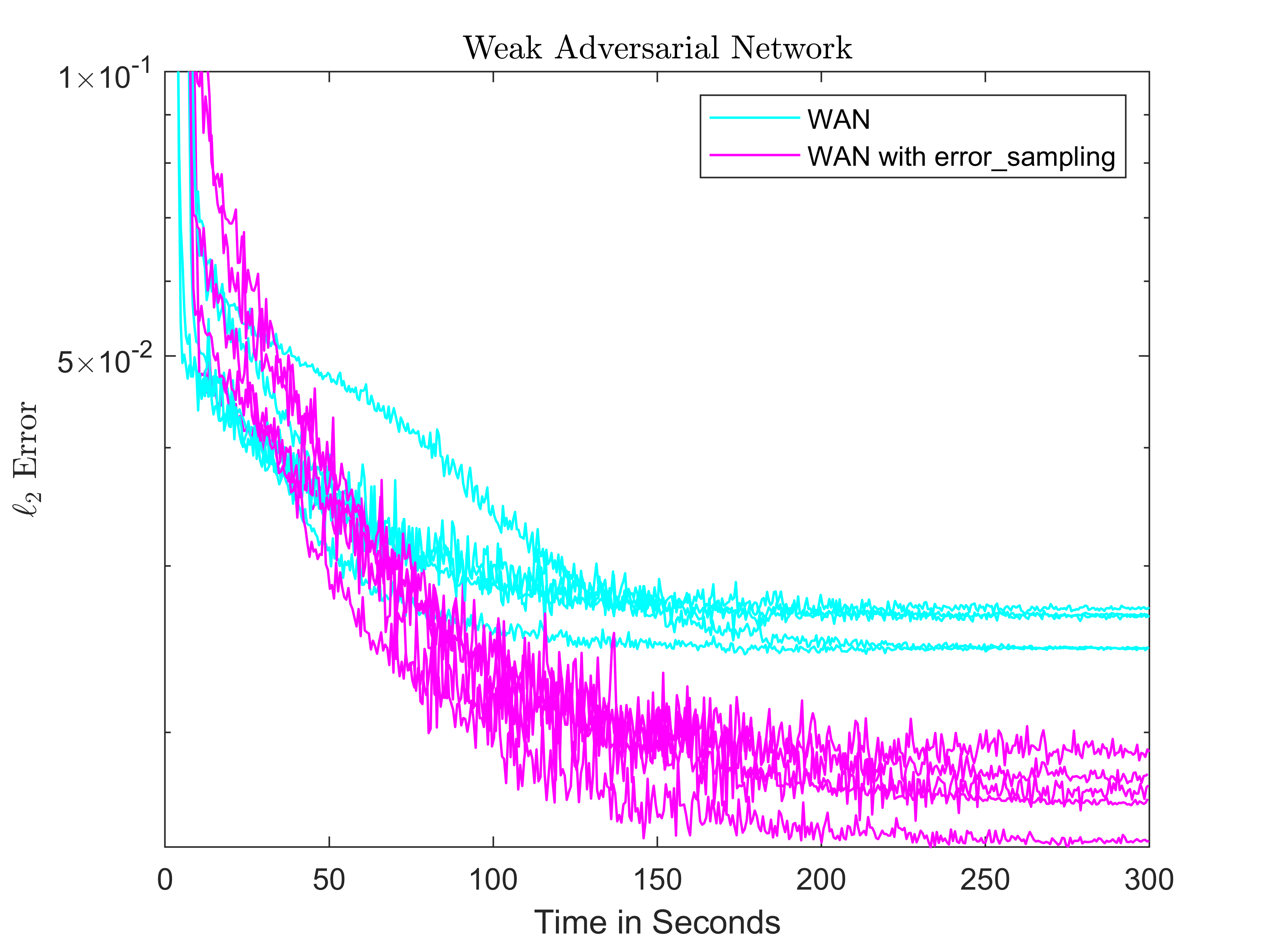}
  \caption*{(c){WAN: $\ell_2$ Error vs Time}}
\endminipage
\end{figure}

Table \ref{table:7} shows statistics of $\ell_2$ errors obtained in 30 trials without and with  adaptive sampling. Figure \ref{fig:5.111} shows the comparison of the five lowest $\ell_2$ error curves obtained with and without adaptive sampling in 30 trials versus time in seconds for DGM, DRM, and WAN respectively. These results clearly show that adaptive sampling is compatible with all three frameworks. The errors obtained with adaptive sampling are lower than those without adaptive sampling. Note that, as observed in previous high-dimensional examples, the advantage of adaptive sampling becomes more significant when the dimension is high.  

\begin{table}[h!]
	\caption{Compatibility tests; $^{\ast}$ denotes training equipped with adaptive sampling}\label{table:7}
	\centering
	\begin{tabular}{lllll}
		\toprule
		Framework & Mean & Standard Deviation & Minimum Value& Coefficient of Variation\\
		\midrule
		DGM & 0.0333 & 0.0064 &  0.0244 & 19.2$\%$  \\
		DGM$^{\ast}$ & 0.0280 & 0.0079 &  0.0155 & 29.2$\%$ \\
		\midrule
		DRM & 0.0273 &0.0097  & 0.0132 & 35.5$\%$     \\
		DRM$^{\ast}$& 0.0255 &0.0093  & 0.0122 & 36.5$\%$\\
		\midrule
		WAN & 0.0329& 0.0062 & 0.0245  & 18.8$\%$ \\
		WAN$^{\ast}$ &0.0282 & 0.0075 &  0.0153 & 26.6$\%$\\
		\bottomrule
	\end{tabular}
\end{table}

\textbf{Acknowledgement} 
W.G. was partially supported by the National Science Foundation under grant DMS-2050133. 
 C. W. was partially supported by the National Science Foundation under awards DMS-2136380 and DMS-2206332. The authors would like to express our deepest appreciation to Dr. Haizhao Yang from  University of Maryland College Park for his invaluable insights, unparalleled support, and constructive advice. Thanks should also go to Dr. Yiqi Gu from University of Hong Kong for helpful discussions.

\bibliographystyle{plainnat}
\bibliography{references}  
\newpage
\begin{appendices}
\section{RAR Sampling}
\begin{algorithm}[H]
\caption{RAR sampling to choose 12000 training points}\label{alg:6.1}
\SetAlgoLined
\setcounter{AlgoLine}{0}
\SetKwInput{KwResult}{Result}
\KwResult{$N= 12000$ points for training}
\SetKwInput{KwRequire}{Require}
\KwRequire{PDE (\ref{eq:2.1}); the current solution net $\phi(\boldsymbol{\theta})$}
 Generate $10000$ uniformly distributed points $\big\{\boldsymbol{x}_i\big\}_{i = 1}^{10000} \subset \Omega$; denote by $\boldsymbol{X}$\;
 Residual$\_$Error$\_$array = $\mathcal{R}^p_{abs}(\boldsymbol{X}) = |\mathcal{D}\phi(\boldsymbol{X}) - f(\boldsymbol{X})|^p$ \;
 Add 2000 points in $\boldsymbol{X}$ with the largest residual errors to $\boldsymbol{X}$ \;
 \Return $\boldsymbol{X}$ for training \;
\end{algorithm}
\noindent The same logic follows for sampling on $\partial \Omega$.

\section{Preliminary Definition}\label{defs}

\begin{definition} Let $\boldsymbol{x}$ be a real random variable and $w(\boldsymbol{x})$ be any real function. The  expectation  of $w(\boldsymbol{x})$ is defined by
$$\mu =\mathbb{E}_{\boldsymbol{\boldsymbol{x}} \in \Omega} \left[w(\boldsymbol{x})\right] =   \int_{\Omega} w(\boldsymbol{x})p(\boldsymbol{x}) \,d\boldsymbol{x}, $$
 where $p(\boldsymbol{x})$ is the probability density function of $\boldsymbol{x}$.
\end{definition}
\begin{definition}
 Let  $\{\boldsymbol{x_i}\}_{i=1}^n \in \Omega \subset \mathbb{R}^d$ be real random variables that follow the distribution $p(\boldsymbol{x})$ and $w(\boldsymbol{x})$ be any function defined on $\Omega$ with density $p(\boldsymbol{x})$. The plain Monte-Carlo estimate of $\mu $, the expectation of $w(\boldsymbol{x})$, is   defined by
$$ \hat{\mu} =  \frac{|\Omega|}{n} \sum_{i=1}^{n} w\left(\boldsymbol{x}_i\right)p\left(\boldsymbol{x}_i\right),   \boldsymbol{x_i} \sim p(\boldsymbol{x}),$$
  where  $|\Omega| = \int_{\Omega}d\boldsymbol{x}$ denotes the volume (length in 1D and area in 2D) of $\Omega$,  $n$ is the number of points sampled for estimating $\mu$.
Note that $\hat{\mu} $ is an unbiased estimate of $\mu$. 
\end{definition}

\section{Numerical Differentiation}\label{num_diff}
\begin{definition}
A real scalar function $\phi(\boldsymbol{x})$ of $d$ variables is a rule that assigns a number $\phi(\boldsymbol{x}) \in \mathbb{R}$ to an array of numbers $\boldsymbol{x} = (x_1, \cdots, x_d)  \in \mathbb{R}^d$.  
\end{definition}

\begin{definition}
Let $\phi(\boldsymbol{x})$ be a real scalar function of $d$ variables defined on $\mathbb{R}^d$. The \textbf{partial derivative} of $\phi(\boldsymbol{x})$ at a point $a \in \mathbb{R}^d$ with respect to $x_i (i = 1, \cdots, d)$ is given by
$$\frac{\partial \phi}{\partial x_i}(a) = \lim_{h_i\to 0} \frac{\phi(a+h\boldsymbol{e_i})-\phi(a)}{h},$$
 where $h \in \mathbb{R}$, $\boldsymbol{e_i} \in \mathbb{R}^d$  is an array of all zeros except 1 for the i-th element.
\end{definition}

\begin{definition}
Let $\phi(\boldsymbol{x})$ be a real scalar function of $d$ variables defined on $\mathbb{R}^d$, $a \in \mathbb{R}^d$, and $h \in \mathbb{R}$. The \textbf{numerical differentiation} estimate of $\frac{\partial \phi}{\partial x_i}(a)$ with respect to $x_i (i = 1, \cdots, d)$ is given by:
$$\frac{\partial \phi}{\partial x_i}(a) \approx \frac{\phi(a+h\boldsymbol{e_i})-\phi(a)}{h}.$$
Note that the truncation error is of order $O(h)$.
\end{definition}

\section{Self-normalized Sampling}\label{self-norm}
\noindent Let $q(\boldsymbol{x})$ be the desired distribution of a continuous variable $\boldsymbol{x}$ over the domain $\Omega$. Suppose that we cannot compute $q(\boldsymbol{x})$ directly but have the unnormalized version $f(\boldsymbol{x})$ such that $$q(\boldsymbol{x}) =  \frac{f(\boldsymbol{x})}{NC}$$
where $NC = \int_{\Omega} f(\boldsymbol{x}) d\boldsymbol{x}$ is the unknown normalizing constant. 
The self-normalized sampling algorithm is:

\begin{algorithm}[H]
\caption{Self-normalized Sampling}\label{alg:d}
\SetAlgoLined
\setcounter{AlgoLine}{0}
\SetKwInput{KwResult}{Result}
\KwResult{$k$ points approximately follow  the target distribution $q(\boldsymbol{x})$}
\SetKwInput{KwRequire}{Require}
\KwRequire{The unnormalized distribution $f(\boldsymbol{x})$}
 Generate an ordered list of $n$ uniformly distributed points $\big\{\boldsymbol{x}_i\big\}_{i = 1}^{n} \subset \Omega$; denote by $\boldsymbol{X}$. \;
 $f(\boldsymbol{X})$ denotes the list $\big\{f(\boldsymbol{x}_i)\big\}_{i = 1}^{n} \subset \Omega$\;
 constant $NC =  \sum_{i=1}^{n} f(\boldsymbol{x}_i)$  \;
 A new discrete probability mass function over $\boldsymbol{X}$ given by $p(\boldsymbol{x}_i) = \frac{f(\boldsymbol{x}_i)}{NC}$\;
 Generate $k$ points following the discrete p.m.f. $p(\boldsymbol{x})$ \;
 \Return $k$ points generated \;
\end{algorithm}
\noindent In the first step, $n$ should be large even if $k$ is small, and when $k$ is large, $n$ can be equal to or even smaller than $k$ as long as $n$ is still large; in the numerical experiment in this work, 12,000 allocation points are needed, so $n$ and $k$ are both set to be 12,000.\\ 

\noindent We will show that this algorithm makes sense. Let $\chi$ be any subset of $\Omega$.\\

\noindent The probability of $\boldsymbol{x} \in \chi$ in the sense of our target distribution $q(\boldsymbol{x})$ is:
\begin{equation}\label{eq:D}
P(\boldsymbol{x} \in \chi) = \int_{\chi} q(\boldsymbol{x}) d\boldsymbol{x} = \int_{\chi} \frac{f(\boldsymbol{x})}{NC} d\boldsymbol{x}  = \frac{\int_{\chi} f(\boldsymbol{x}) d\boldsymbol{x}}{\int_{\Omega} f(\boldsymbol{x}) d\boldsymbol{x}} \approx \frac{\int_{\chi}d\boldsymbol{x}\cdot \frac{1}{m} \sum_{j=1}^{m} f(\boldsymbol{x}_j)}{\int_{\Omega}d\boldsymbol{x}\cdot \frac{1}{n} \sum_{i=1}^{n} f(\boldsymbol{x}_i)} = \frac{n}{m} \cdot \frac{\int_{\chi}d\boldsymbol{x}}{\int_{\Omega}d\boldsymbol{x}} \cdot \frac{\sum_{j=1}^{m} f(\boldsymbol{x}_j)}{\sum_{i=1}^{n} f(\boldsymbol{x}_i)}.
\end{equation}

\noindent Note that, this is the Monte Carlo integration, $\boldsymbol{x}_j$ and $\boldsymbol{x}_i$ are $m$ and $n$ random points in $\chi$ and $\Omega$ respectively. This approximation is governed by the Law of Large Numbers; therefore, $m$ and $n$ are expected to be large.\\

\noindent Now, suppose out of $n$ points generated in Algorithm \ref{alg:d}, there are $m$ points in $\chi$. When $m$ and $n$ are large, by the Law of Large Numbers, $\frac{m}{n} \approx \frac{\int_{\chi}d\boldsymbol{x}}{\int_{\Omega}d\boldsymbol{x}}$. Hence, we can apply these $n$ and $m$ as random points to Equation \ref{eq:D}, the above probability becomes:
\begin{equation}
P(\boldsymbol{x} \in \chi) \approx \frac{\sum_{j=1}^{m} f(\boldsymbol{x}_j)}{\sum_{i=1}^{n} f(\boldsymbol{x}_i)},
\end{equation}
  and this is the probability of $\boldsymbol{x} \in \chi$ in the sense of the discrete p.m.f. $p(\boldsymbol{x})$ in Algorithm \ref{alg:d}.
  
One advantage of this algorithm is that it generates points in a completely parallel fashion whereas in the Metropolis-Hastings algorithm, each draw is based on the previous draw; thus, the MH algorithm cannot be parallelized.
\end{appendices}

\end{document}